%% file: main.tex
\documentclass{article}

\usepackage{arxiv}

\usepackage[utf8]{inputenc} % allow utf-8 input
\usepackage[T1]{fontenc}    % use 8-bit T1 fonts
\usepackage{hyperref}       % hyperlinks
\usepackage{url}            % simple URL typesetting
\usepackage{booktabs}       % professional-quality tables
\usepackage{amsfonts}       % blackboard math symbols
\usepackage{nicefrac}       % compact symbols for 1/2, etc.
\usepackage{lipsum}		    % Can be removed after putting your text content
\usepackage{graphicx}
\usepackage[square,sort,comma,numbers]{natbib}
\usepackage{doi}
\usepackage{float}
\usepackage{caption}
\usepackage{subcaption}
\usepackage{xcolor}
\usepackage{amsmath}
\usepackage{amssymb}
\usepackage{svg}
\usepackage{enumitem}
\usepackage{multirow}

\title{An unstructured second-order subgrid method for the shallow water equations}

\author{ \href{https://orcid.org/0009-0003-4603-5073}{\includegraphics[scale=0.06]{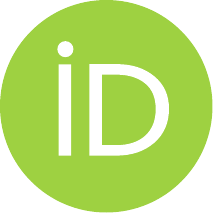}\hspace{1mm}Max Ebstrup Bitsch} \\
Department of Applied Mathematics and Computer Science\\
Technical University of Denmark\\
Kongens Lyngby, 2800, Denmark, \\
and DHI A/S\\
Hørsholm, 2970, Denmark \\
\texttt{mebi@dhigroup.com} \\
\And
    \href{https://orcid.org/0000-0003-2185-1165}{\includegraphics[scale=0.06]{Figures_pdf/orcid-eps-converted-to.pdf}\hspace{1mm}Irene Torpe Heilmann} \\
    DHI A/S\\
    Hørsholm, 2970, Denmark\\
\And
    \href{https://orcid.org/0000-0001-8626-1575}{\includegraphics[scale=0.06]{Figures_pdf/orcid-eps-converted-to.pdf}\hspace{1mm}Allan Peter Engsig-Karup} \\
    Department of Applied Mathematics and Computer Science\\
    Technical University of Denmark\\
    Kongens Lyngby, 2800, Denmark \\
\And
    Ole Rene Sørensen \\
    DHI A/S\\
    Hørsholm, 2970, Denmark\\
\And
    \href{https://orcid.org/0000-0002-4331-9244}{\includegraphics[scale=0.06]{Figures_pdf/orcid-eps-converted-to.pdf}\hspace{1mm}Jesper Grooss} \\
    DHI A/S\\
    Hørsholm, 2970, Denmark\\
}
\renewcommand{\headeright}{Bitsch et al. (2026)}
\renewcommand{\shorttitle}{An \textit{arXiv} Preprint}

\hypersetup{
pdftitle={},
pdfsubject={},
pdfauthor={},
pdfkeywords={},
}

\begin{document}
\maketitle

\begin{abstract}
	\input{Texts/Abstract.tex}
\end{abstract}

% keywords can be removed

% Input your different sections here:
\input{Texts/Introduction}
\input{Texts/Mathematical_model}
\input{Texts/Mesh}
\input{Texts/Numerical_discretization}

\input{Texts/Results_verification}
\input{Texts/Conclusion_and_more}

\bibliographystyle{SageH}

% When working on the draft use the following:
\bibliography{references.bib}  

% When submitting to arXiv (example) - remember to select all .bib etc. files.:



% \end{thebibliography}

\appendix
\input{Texts/Appendix}

\end{document}

%% file: Texts/Abstract.tex
% NOTE: MAXIMUM 200 words!
We propose a new unstructured numerical subgrid method for solving the shallow water equations using a finite volume method with enhanced bathymetry resolution. 
The method employs an unstructured triangular mesh with support for triangulation of elements to form finer subgrids locally. The bathymetry is represented on the fine mesh, allowing the incorporation of small-scale features, while velocities are defined on the coarse mesh. The governing equations are solved numerically on the coarse mesh, making the method computationally cheaper compared to a traditional fine mesh computation. To accurately represent the velocities, we employ a second-order accurate WENO method and for temporal integration an explicit second-order accurate Runge-Kutta method. Furthermore, we present a novel subgrid face value reconstruction that accounts for partially wet cells, where only some of the subgrid cells are wet. 
Together with a newly developed gravity source term discretization, we demonstrate that the scheme is well-balanced on the subgrid level. 
Finally, the subgrid method is validated on several test cases to show: i) that the scheme is well-balanced, ii) confirm the order of the accuracy of the scheme, and iii) demonstrate the capability of handling moving flood and dry boundaries. We also highlight when it is beneficial to increase the number of subgrid cells to improve the results of the finite volume method.

%% file: Texts/Introduction.tex
\section{Introduction} \label{sec:Intro}
Flooding poses a major threat in today's society as it can cause immense amounts of damage to critical infrastructure. This is an increasing problem as floods have become more frequent over the last decades \cite{MADSEN2014, AGONAFIR2023}. Flood management is key in safety assessments aimed at reducing damage, where early risk management is especially important. Here, flood models play a crucial role in predicting the impact of floods and forecasting potential damage. With such information, early warning systems and other infrastructure can be built to reduce or even prevent damage from floods.   

Many flood models are based on the two-dimensional shallow water equations (SWE). The SWE are depth-integrated Navier-Stokes equations that -- in shallow regions -- lead to a solid compromise between model accuracy and computational efficiency. This is why they are very popular in commercial flooding applications. 
SWE are often discretized using the Finite Volume method (FV) due to the hyperbolic nature of the governing equations \cite{Kurganov2018,Delgrosso2024}. FV methods are inherently mass conserving and can accurately propagate shocks when implemented correctly, which are important properties for water flow modeling. For flood modeling, the treatment of wetting and drying is particularly important but non-trivial to incorporate, and special care must be taken. The challenge lies in accurately predicting the wetting and drying front and determining the correct velocity as the water depth approaches zero while maintaining model stability. Here, we provide a brief overview of representative approaches \cite{Begnudelli2006,Song2011,Hou2013,Cea2015,Xia2017,CaviedesVoullieme2023,Gallardo2007}.

Inherent in these numerical models is the fact that an increase in accuracy corresponds to an increase in computational requirements. 
This has become a problem, as the desire to resolve fine features such as sidewalk curbs and other small obstacles, is now possible as there are more satellite data available than ever before.
To combat long run times, the use of modern computer architecture, such as multiple CPUs or GPUs, can be employed \cite{BrodtkorbEtAl2011,MORALESHERNANDEZ2021,CaviedesVoullieme2023}. 
An alternative to this "brute force" approach of scaling hardware resources for improved run times is instead an algorithmic approach, where the numerical method is tailored to the problem. Traditionally, unstructured meshes have been used because they are boundary-fitted, avoiding a staircase representation of domain boundaries, and they can be refined locally in regions with complex features such as small streams or urban buildings. Local refinement greatly reduces the number of computational cells compared with a structured mesh, where every cell has the same size. However, for large areas with many features, the scale of the computations can still become too large for practical run times that meet time budget constraints required for simulations and analysis.
Another candidate for an algorithmic approach is the so-called {\em subgrid methods}, where a computational cell has multiple bathymetry values. This accounts for local bathymetry variations that would otherwise be obscured.
Inclusion of additional bathymetry allows for local flooding and drying on subgrid level. This allows a computational cell to be partially wet, in contrary to a traditional method where it is either wet or dry. Ultimately, this means that mesh-fitting can be accounted for on a subgrid level.   

One type of subgrid method are the porosity models, mostly used for urban areas. These methods consider the wet part of a cell, based on a wet/dry fraction. It was initially proposed by \cite{Defina1994} with more comprehensive details given later \cite{Defina2000}. It has since been extended to a wide range of applications, some of which even include depth-dependencies \cite{Guinot2006,Sanders2008,OZGEN2016,GUINOT2018}.

A more direct depth-dependent method using high-resolution bathymetry data for the SWE, proposed by \cite{Casulli2009}, is referred to as the Casulli subgrid method. The method is a semi-implicit scheme and is based on orthogonal meshes, making the mesh flexible but not completely unstructured. The mesh is combined with an underlying subgrid mesh to provide a more accurate determination of the water volume and the wet/dry front. It was later extended to three dimensions \cite{Casulli2011}. The 2D method has inspired many different articles in this area \cite{Volp2013,Sehili2014,KENNEDY2019,Zhang2021,STELLING2022}. 
To overcome the limitation of having coarse velocity resolution with restricted meshing it is possible to use quadtrees \cite{Stelling2012}.

More recently, the Casulli subgrid method has been extended to -- what the author calls -- grid cloning \cite{Casulli2019}, where multiple rivers and streams can flow into the same element face and remain separated. This allows for very coarse meshes and very low run times for certain types of applications. 

A depth-dependent subgrid method that is explicit and solves the Riemann problem was proposed by \cite{Sanders2019}. This method has the underlying subgrid mesh, but is restricted to a Cartesian mesh. It uses a second-order reconstruction of the free surface to enhance the results. A similar methodology, but with face averaging and Taylor expansions, was introduced in \cite{Shamkhalchian2021}. This method has been extended to second-order in space \cite{Shamkhalchian2023}, where they highlight the benefits with respect to computational efficiency.

\subsection{Paper contribution}
The primary objective of this work is to propose a new cell-centered FV subgrid method that utilizes unstructured triangular meshes. This approach enables features smaller than a mesh element to be represented through the subgrid while ensuring that boundaries are correctly fitted and velocities can be locally refined to accurately capture the velocity field.  
The work includes, a novel subgrid face value reconstruction to account for local subgrid variation including flood and dry. This reconstruction and a newly developed gravity discretization makes the method well-balanced on subgrid level. The scheme is implemented using the weighted essentially non-oscillatory (WENO) method making it second-order accurate, enabling it to represent velocities more accurately.  
The scheme is verified primarily through test cases from \cite{Delestre2013} that include numerous analytical solutions. This highlights that important numerical and physical features are preserved.    

\subsection{Paper outline}
Following the introduction in Section \ref{sec:Intro}, the governing model equations are presented in Section \ref{sec:Gov_Eq}, including choices of friction and gravitational representation. Section \ref{sec:Subgrid_Mesh} describes the domain, the subgrid meshing with its coarse and fine mesh, its notation, and how the variables are represented. Subsequently, Section \ref{sec:Num_Meth} outlines the numerical spatio-temporal discretization, highlighting the second-order accurate (WENO) method with a focus on flood and dry regions. This section also includes a description of the temporal time-stepping approach, the face value reconstruction, the numerical flux, and the discretization of gravity and friction. Results are presented in Section \ref{sec:Model_Verfication}, with an emphasis on cases that have analytic solutions, as well as two cases that highlights the subgrid methods' ability to handle small-scale features. Finally, the paper is concluded with a summary of the main findings in Section \ref{sec:Conclusion}.    

%% file: Texts/Mathematical_model.tex
\section{Governing equations} \label{sec:Gov_Eq}
\begin{figure}[htb]
    \centering
    \includegraphics[width=0.55\linewidth]{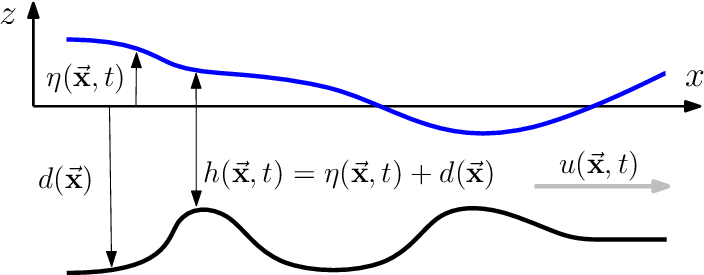}
    \caption{Illustration of the computational domain and notation for the shallow water equations.}
    \label{fig:SWE_illu}
\end{figure}

The two-dimensional SWE acting on the domain, $\Omega \subset \mathbb{R} ^2$, with $\vec{\mathbf{x}} = (x,y) \in \Omega$ and a time interval with coordinate $t \in \mathbb{R}_0^+$. The SWE are a system of first-order hyperbolic differential equations given by
\begin{equation}
     \frac{\partial \mathbf{U}}{\partial t} +   \nabla \cdot \mathbf{F}(\mathbf{U}) =  \mathbf{S}(\mathbf{U}), \label{eq:gov_eq}
\end{equation}
where $\mathbf{U} = (h,hu,hv)^T$ are the conserved variables, with $h = h(\vec{\mathbf{x}},t)$ being the water depth, $u = u(\vec{\mathbf{x}},t)$ and $v = v(\vec{\mathbf{x}},t)$ being the velocities in the $x$- and $y$-direction, respectively. The water depth connects the free surface and the bathymetry as $h = \eta + d$, where $\eta = \eta(\vec{\mathbf{x}},t)$ is the free surface elevation and $d = d(\vec{\mathbf{x}})$ is the time-invariant bathymetry. Note that the bathymetry is positive downwards, as shown in Figure \ref{fig:SWE_illu}, which illustrates a cross-section along the $x$-axis.    

\noindent The flux, $\mathbf{F}(\cdot)$, is given by the advection term and part of the gravity term as
\begin{equation}
    \mathbf{F}(\mathbf{U}) = \begin{pmatrix}
        hu & hv \\
        hu^2 + \frac{1}{2} g (\eta^2 + 2\eta d)  & huv \\
        hvu & hv^2 + \frac{1}{2} g (\eta^2 + 2\eta d)
    \end{pmatrix}, \label{eq:flux}
\end{equation}
where $g=9.81$ ms\textsuperscript{-2} is the gravitational acceleration assumed to be constant. The source term is given as the sum of the external forces, gravity and friction, as
\begin{equation}
        \mathbf{S}(\mathbf{U}) = \mathbf{S}_g(\mathbf{U}) + \mathbf{S}_f(\mathbf{U}) = \begin{pmatrix}
        0 \\
        g\eta \frac{\partial d}{\partial x} \\
        g\eta \frac{\partial d}{\partial y}
    \end{pmatrix} +
    \begin{pmatrix}
        0 \\
        - c_{fb} |\vec{\mathbf{u}}| u \\
        - c_{fb} |\vec{\mathbf{u}}| v 
    \end{pmatrix}, \label{eq:source}
\end{equation}
where $c_{fb}=g/(M^2 h^{1/3})$ is the friction coefficient given by the Mannings coefficient, $M$, here being larger than one ($M \geq 1$), and $\vec{\mathbf{u}} = (u,v)^T$ is the velocity vector with $|\vec{\mathbf{u}}| = \sqrt{u^2+v^2}$ being the magnitude. The gravity term follows the reformulation of \cite{Liang2009}, which helps balance the numerical flux and the source term. The expression is divided such that one part appears in the flux term and the other in the source term.

%% file: Texts/Mesh.tex
\section{Domain, Variables, and Subgrid mesh} \label{sec:Subgrid_Mesh}
The domain is partitioned into two meshes: an unstructured coarse mesh, called the coarse mesh, and an underlying fine mesh, called the subgrid mesh. 

Let the true spatial domain $\Omega_{true} \subset \mathbb{R}^2$ be approximated by the computational domain $\Omega \simeq \Omega_{true}$ divided into $N_C$ unstructured non-overlapping triangular elements $\tilde{T}_m, ~m = 1,2,..., N_c$, such that
 \begin{equation}
    \Omega = \bigcup_{m=1}^{N_c} \tilde{T}_m,
 \end{equation}
called the coarse mesh denoted by the tilde.
Note that if the true domain is a polygon it can be covered exactly $\Omega = \Omega_{true}$. Now we partition each element into $N_m^{sg}$ subgrid elements $T_{m,k}$ which are also non-overlapping triangular elements such that $\tilde{T}_m = \bigcup_{k=1}^{N^{sg}_m} T_{m,k}$, where $k$ is local numbering. From this, it follows that the entire domain can be represented by subgrid elements such that $\Omega = \bigcup_{l=1}^{N_{sg}} T_l$, where $N_{sg}$ is the total number of subgrid elements in the entire domain and $l$ is the global numbering of subgrid elements. This is called the global subgrid mesh, and a mapping from the local to the global subgrid mesh is given by
\begin{equation}
    l = \ell(m,k), \quad k = 1,2,...,N^{sg}_m.
\end{equation}
From this, the domain can also be represented as
\begin{equation}
    \Omega = \bigcup_{m=1}^{N_c}\bigcup_{k=1}^{N^{sg}_m}T_{\ell(m,k)}.
\end{equation}

Each triangular element has three non-curved edges/faces, which together constitute the boundary of the element 
\begin{equation}
    \partial \tilde{T}_m = \bigcup_{i=1}^3\tilde{E}_{m,i},
\end{equation}   
where $\tilde{E}$ is the edges ordered counterclockwise. Similarly, for the subgrid elements we have that 
\begin{equation}
    \partial T_{\ell(m,k)} = \bigcup_{j=1}^3 E_{\ell(m,k),j}.
\end{equation}
For every coarse element edge, there exists a set of subgrid edges given by their indices $k$ and $j$ such that they overlap with a coarse face 
\begin{equation}
    \mathcal{E}_{m,i} = \{(k,j) | \tilde{E}_{m,i} \cap E_{\ell(m,k),j}, k = 1,2,...,N_m^{sg} ~\text{and} ~ j= 1,2,3 \},
\end{equation}
and where the union spans the $i$'th coarse edge 
$\tilde{E}_{m,i} = \bigcup_{(k,j) \in \mathcal{E}_{m,i}} E_{\ell(m,k),j}$. In this work, we have chosen to divide each coarse face into the same number of subgrid faces $n_m^{sg}$, with equal length. The set $\mathcal{E}_{m,i}$ is then not unique as the edges are shared by two elements -- element $m$ and its neighbor across the $i$'th edge $m_i$ -- with similar subgrid face division. We get that $\{E_{\ell(m,k),j}|(k,j)\in \mathcal{E}_{m,i} \} = E_{\ell(m_i,k),j}|(k,j)\in \mathcal{E}_{m_i,i^{\prime}} \}$, where $i^\prime \in \{1,2,3\}$. 
% $\mathcal{E}_{m,i} = \mathcal{E}_{m_i,i^\prime}$, where $i^\prime \in \{1,2,3\}$.  

\begin{figure}[htb]
    \centering
    \includegraphics[width=0.65\linewidth]{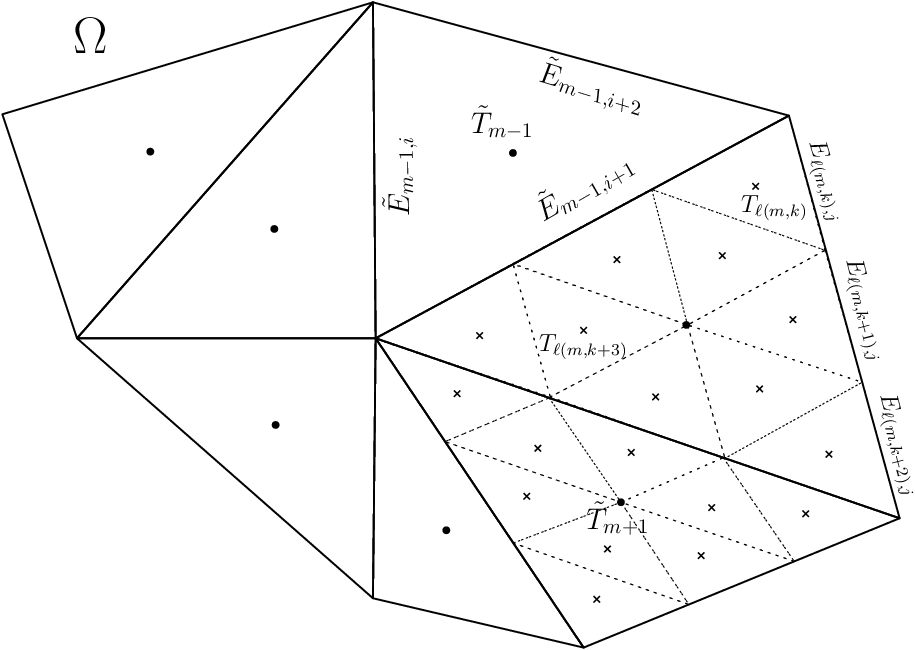}
    \caption{Illustration of the computational domain $\Omega$ and the mesh notation, where the domain is partitioned into $N_C$ coarse unstructured triangles $\tilde{T}_m$ (solid black lines), with some being divided into $N_m^{sg}$ subgrid elements $T_{\ell(m,k)}$ (dashed lines).}
    \label{fig:Mesh_con}
\end{figure}

An illustration of a cell-centered, unstructured triangular mesh is shown in Figure \ref{fig:Mesh_con}, along with the corresponding notation. In the figure, only two cells are illustrated with subgrid resolution -- this is purely for illustrative purposes -- where a typical mesh would contain subgrid elements in all coarse elements. The coarse cells are divided into $N_{m}^{sg}=9$ subgrid cells with $n_m^{sg}=3$, where the element $\tilde{T}_m$ also shows the subgrid faces included in the set $\mathcal{E}_{m,i}$. 

\subsection{Variable representation}
As for cell-centered FV method, the variables are represented as integral averaged quantities at the cell centers; hence, the coarse mesh variables are
\begin{equation}
    \bar{\mathbf{U}}_m = \frac{1}{|\tilde{T}_m|}\int_{\tilde{T}_{m}} \mathbf{U}~\text{d} \vec{\mathbf{x}},
\end{equation}
where $|\tilde{T}|$ is the coarse cell area. In order to benefit from the additional subgrid resolution, the discrete integral average becomes
\begin{equation}
    \bar{\mathbf{U}}_m = \frac{1}{|\tilde{T}_m|}\sum_{k=1}^{n_m^{sg}} \int_{T_{\ell(m,k)}} \mathbf{U}~\text{d} \vec{\mathbf{x}} = \frac{1}{|\tilde{T}_m|} \sum_{k=1}^{N_m^{sg}} |T_{\ell(m,k)}| \mathbf{U}_{\ell(m,k)}. \label{eq:Int_Avg}
\end{equation}
Here, we have chosen the free surface elevation and velocities as reconstruction variables, which are based on the coarse mesh. These are stored in the vector $\mathbf{Q}$, which are determined at the cell-centers as $\bar{\mathbf{Q}}_m = (\bar{\eta}_m,\bar{hu}_m/\bar{h}_m,\bar{hv}_m/\bar{h}_m)^T$, with $\bar{\eta}_m$ being defined in Section \ref{sec:Discrete_Mass}.

The subgrid reconstruction is based on polynomials, and the subgrid cell-center values are 
\begin{subequations}
    \begin{align}
        \eta_{\ell(m,k)} =  \mathbf{w}_m^{\eta}(\vec{\mathbf{x}}_{\ell(m,k)}), \\ u_{\ell(m,k)} =  \mathbf{w}_m^{u}(\vec{\mathbf{x}}_{\ell(m,k)}), \\  v_{\ell(m,k)} =  \mathbf{w}_m^{v}(\vec{\mathbf{x}}_{\ell(m,k)}),
    \end{align}
\end{subequations}
where $\mathbf{w}_m^{\cdot}$ is a polynomial on the $m$'th cell based on the super-scripted variable and $\vec{\mathbf{x}}_{\ell(m,k)}$ are coordinates of the subgrid cell centers. For this work, the polynomials are limited to being either zeroth- or first-order ($P_n=0$ or $P_n=1$), which leads to first- and second-order accuracy, respectively. In the zeroth-order case, the polynomial is constant and equal to the coarse cell center value everywhere in the element. The first-order polynomial is determined with the WENO method described in Section \ref{sec:WENO}.
The bathymetry values are stored at the subgrid level and are always piecewise constant, also when using first-order polynomials.

From this, we obtain subgrid variation of the water depth and the flow discharges. The water depth is defined on the subgrid level to be non-negative
\begin{equation}
    h_{\ell(m,k)} = \max(0, \mathbf{w}_m^{\eta}(\vec{\mathbf{x}}_{\ell(m,k)}) + d_{\ell(m,k)}), \label{eq:hsg}
\end{equation}
and the discharges are computed as 
\begin{subequations}
\begin{align}
    hu_{\ell(m,k)} = h_{\ell(m,k)}\mathbf{w}_m^{u}(\vec{\mathbf{x}}_{\ell(m,k)}), \\ hv_{\ell(m,k)} = h_{\ell(m,k)}\mathbf{w}_m^{v}(\vec{\mathbf{x}}_{\ell(m,k)}),
\end{align} 
\end{subequations}
where we have chosen to represent the velocities on the coarse mesh, but this could be replaced by the discharges. 

\begin{figure}[htb]
    \centering
    \subfloat[]{
    \includegraphics[width=30mm]{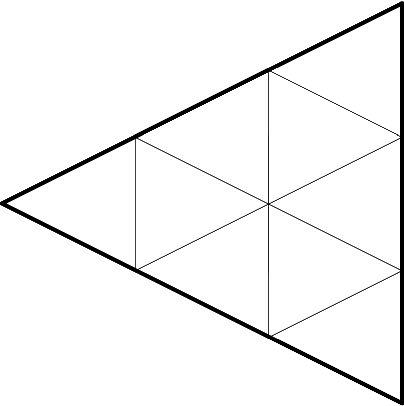}
  % \includesvg[width=30mm]{Figures/DC.svg}
}
    \hspace{20mm}
    \subfloat[]{
    \includegraphics[width=30mm]{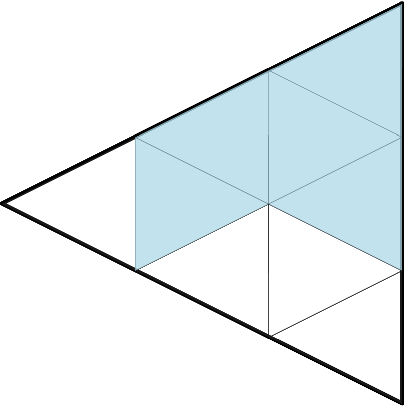}
  % \includesvg[width=30mm]{Figures/PC.svg}
}
    \hspace{20mm}
    \subfloat[]{
    \includegraphics[width=30mm]{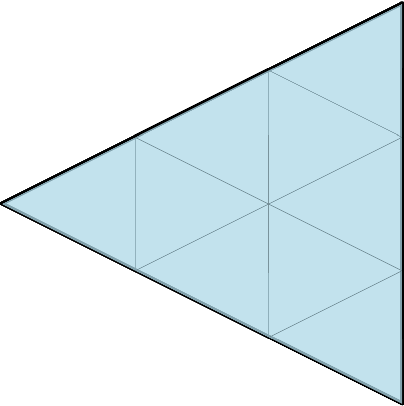}
  % \includesvg[width=30mm]{Figures/WC.svg}
}
    \caption{Illustrations of a coarse cell with $N^{sg}_m = 9$ subgrid cells and the three different cell states, where (a) is dry, (b) is partially wet, and (c) is wet. 
    } \label{fig:Cell_perm}
    
\end{figure}

The non-negative subgrid water depth means that a coarse cell can now contain both wet and dry areas, where a wet subgrid cell is given by $h_{\ell(m,k)} > 0$ and a dry is given by $h_{\ell(m,k)} = 0$. This leads to a new cell state, which we refer to as partially wet. We will distinguish between three different states: \textit{a}) a dry cell, where every subgrid cell is dry, \textit{b}) a partially wet cell, where some subgrid cells are wet but not all, and \textit{c}) a wet cell, where every subgrid cell is wet. Illustrations of all states can be seen in Figure \ref{fig:Cell_perm}, where three coarse elements are shown, with white subgrid cells being dry and blue being wet.    

In order to avoid artificially large velocities in very shallow regions, we use a velocity cut-off approach. This approach may limit the transport of water between cells, but the scheme remains mass conserving. The cut-off approach is chosen due to its simplicity. For a more comprehensive discussion, see \cite{Kurganov2007}.
Here, velocities are set to zero if the coarse water depth is beneath a threshold
\begin{equation}
    \vec{\mathbf{u}}_m(\vec{\mathbf{x}}) = \begin{cases}
    0   \quad &\text{if}~ \tilde{h}_m < \epsilon_{dry} \\
    \tilde{hu}_m/\tilde{h}_m \quad &\text{elsewhere}
    \end{cases},
\end{equation} 
where $\epsilon_{dry}$ is the dry tolerance. 
In cases where the coarse water depth is negative, water from neighboring cells are redistributed to preserve physical water levels.   

%% file: Texts/Numerical_discretization.tex
\section{Numerical method} \label{sec:Num_Meth}

We integrate the governing equations \eqref{eq:gov_eq} over each element and apply the Gauss divergence theorem to get a FV form 
\begin{equation}  
    \int_{\tilde{T}_{m}}\frac{\partial \mathbf{U}}{\partial t} ~\text{d}\vec{\mathbf{x}} +  \int_{\partial \tilde{T}_{m}}  \mathbf{F}(\mathbf{U})\cdot \mathbf{n}~\text{d}\vec{\mathbf{s}} =  \int_{\tilde{T}_{m}} \mathbf{S}(\mathbf{U})~\text{d}\vec{\mathbf{x}},
\end{equation}
where $\vec{\mathbf{s}} = \vec{\mathbf{x}}\in\partial \tilde{T}_m$ are the coordinates on the edges of the \textit{m}'th element and $\mathbf{n} = (n^x,n^y)$ is the outward pointing normal. We apply a one-point quadrature rule to obtain second-order accuracy of the integrals
\begin{equation}  
    \frac{\partial \bar{\mathbf{U}}_m}{\partial t} = - \frac{1}{|\tilde{T}_{m}|} \sum_{i=1}^3\sum_{(k,j)\in\mathcal{E}_{m,i}}  \mathcal{F}(\bar{\mathbf{U}}_{L},\bar{\mathbf{U}}_{R},\mathbf{n}_{m,i})|E_{\ell(m,k),j}| +  \mathcal{S}(\bar{\mathbf{U}}_{m}),
\end{equation}
where $|E|$ is the subgrid edge length, $\mathcal{S}(\cdot)$ is the numerical source term, $\mathcal{F}(\cdot)$ is the numerical flux function that is either the solution to the Riemann problem or a wall flux. The reconstructed face values from the $k$'th subgrid element and the neighbor subgrid element across the $j$'th edge are denoted by the subscript $L$ and $R$, respectively. Here, we have that $L$ always refers to the element in question and $R$ its neighboring element at a cell face. The reconstruction procedure used to determine the face values is described in Section \ref{sec:ValRec}.
For simplicity the flux term is represented as
\begin{equation}
    \hat{\mathcal{F}}(\bar{\mathbf{U}}_m) = - \frac{1}{|\tilde{T}_{m}|} \sum_{i=1}^3\sum_{(k,j)\in\mathcal{E}_{m,i}}  \mathcal{F}(\bar{\mathbf{U}}_{L},\bar{\mathbf{U}}_{R},\mathbf{n}_{m,i})|E_{\ell(m,k),j}|.
\end{equation}

\subsection{The weighted essentially non-oscillating method} \label{sec:WENO}
Second-order accuracy is obtained by using the WENO method. The method we use here is based on the work \cite{Dumbser2007}, where polynomials are reconstructed instead of the point values. The polynomial reconstruction enables efficient evaluation of values at subgrid cell centers and faces. Since the polynomial spans the entire coarse cell, the reconstruction is computed only once per coarse cell, whereas a point-wise reconstruction would require multiple evaluations for each subgrid cell.

\begin{figure}[htb]
    \centering
    \subfloat[]{
    \includegraphics[width=45mm]{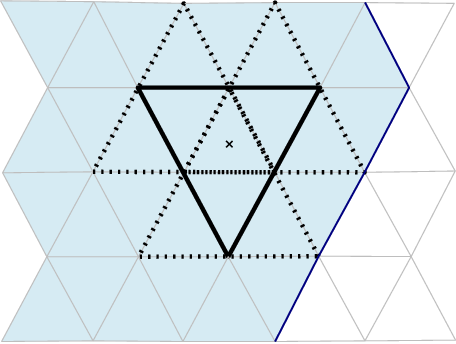}
    \label{fig:Stencils_a}
}
    \hspace{0mm}
    \subfloat[]{
    \includegraphics[width=45mm]{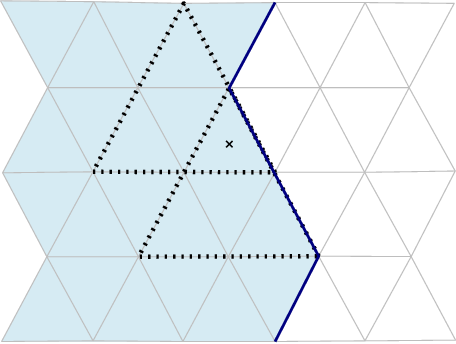}
    \label{fig:Stencils_b}
}
    \hspace{0mm}
    \subfloat[]{
    \includegraphics[width=45mm]{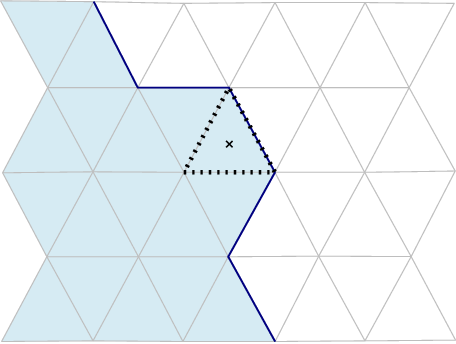} \label{fig:Stencils_c}
}
    \caption{Illustrations of stencils near the flood and dry boundary, where (a) is the complete set of all four stencils, (b) is reduced to two wet sector stencils, and (c) is with zero wet stencils, reducing the scheme to first-order.} \label{fig:Stencils}  
\end{figure}

We restrict it to first-order polynomials (linear planes), which we represent as
\begin{equation}
    \mathcal{P}_m(\vec{\mathbf{x}}) = \bar{\mathbf{Q}}_m + \nabla \mathbf{Q}_m \cdot  (\vec{\mathbf{x}} - \vec{\mathbf{x}}_m),
\end{equation}
where $\nabla \mathbf{Q} = (\partial\mathbf{Q}/\partial x,\partial\mathbf{Q}/\partial y)$ are gradient coefficients of the plane, and $\vec{\mathbf{x}}_m$ are coordinates of the coarse cell centers. The coefficients are determined based on integral values in a set of neighboring cells, referred to as the stencil. In order to get a non-oscillating polynomial, we construct four different polynomials based on four different stencils: the center stencil and three sector stencils. The center stencil is the optimal one when the solution is smooth and is based on the three direct neighbors. The sector stencils are based on cells in a region spanned by the vertices and can be considered as an "up-winding" approach. Each of the four stencils consists of four elements: the target element, on which the polynomial is constructed, and its three neighboring elements. The four stencils can be seen in Figure \ref{fig:Stencils_a}, where the polynomial is constructed in the cell marked with the cross. 
The coefficients are determined by requiring the polynomial to satisfy integral conservation in all cells in the stencil
\begin{equation}
    \int_{\tilde{T}_{\xi(\kappa)}}\mathcal{P}_m(\vec{\mathbf{x}}) ~\text{d} \vec{\mathbf{x}} = \bar{\mathbf{Q}}_{\xi(\kappa)},
\end{equation}
where $1 \leq \kappa \leq 4$ is the local index for the elements in the stencil and $m=\xi(\kappa)$ is the mapping from local to global index. This leads to an overdetermined system, which is solved in a least-squares sense; see \cite{Dumbser2007} for more details. The final WENO polynomial is constructed as a weighted linear combination of the four different stencils 
\begin{equation}
    \mathbf{w}_m(\vec{\mathbf{x}}) = \sum_{s=0}^{3} \omega_m^s \mathcal{P}_m^s(\vec{\mathbf{x}}),
\end{equation}
where $\omega_s$ are nonlinear normalized weights defined as
\begin{equation}
    \omega_m^s = \frac{\tilde{\omega}_m^s}{\sum_{s=0}^{3}\tilde{\omega}_m^s}, \quad \text{where} \quad \tilde{\omega}_m^s = \frac{\lambda_s}{(\sigma_m^s + \epsilon)^p},
\end{equation}
with $\sigma$ being an smoothness indicator, $\lambda $ the linear weights, $\epsilon$ a small parameter to avoid division by zero, and $p$ is a parameter to strengthen the effects of large smoothness indicator values. The linear weights are determined in such a way that they favor the center polynomial ($s=0$) to ensure optimal accuracy. Here we use the values from \cite{Dumbser2017}
\begin{equation}
    \lambda_s = \begin{cases} \lambda_c\quad  &s  = 0 \\ 
        1 \quad &\text{otherwise} 
    \end{cases}~,
\end{equation}
with $\lambda_c = 10^5$. For the remaining parameters we have $\epsilon = 10^{-14}$ and $p=4$. The smoothness indicator is computed in accordance with \cite{JIANG1996,HU1999}
\begin{equation}
    \sigma_m^s = \sum_{1 \leq |\alpha|} \int_{\tilde{T}_m} \left( \frac{\partial^{|\alpha|} \mathcal{P}_m^s(\vec{\mathbf{x}}) }{\partial x^{\alpha_1}\partial y^{\alpha_2}} \right)^2 |\tilde{T}_m|^{|\alpha| - 1} ~\text{d} \vec{\mathbf{x}},
\end{equation}
where $\alpha$ is a multi-index describing the order of derivative, and $|\alpha|$ is the sum of the indices. Fortunately, this reduces to a more approachable expression in the case of a first-order polynomial
\begin{equation}
    \sigma_m^s = |\tilde{T}_m|\left(\left(\frac{\partial \mathbf{Q}_m}{\partial x}\right)^2 + \left(\frac{\partial \mathbf{Q}_m}{\partial y}\right)^2\right). 
\end{equation}
It should be noted that the smoothness indicator and the nonlinear weights are computed for each polynomial for each reconstruction variable.

\subsubsection{Stencils near dry cells}

In the presence of dry cells, the polynomials become non-physical as the cells have no meaningful value. Here, we employ a simple yet effective approach, utilizing only completely wet stencils. A completely wet stencil contains wet and partially wet cells but no dry ones. This means that if just a single coarse cell in a stencil is dry, the entire stencil is discarded. In the case where all stencils are discarded, the cell value is used, making it locally first-order accurate. An example of a reduced stencil set is shown in Figure \ref{fig:Stencils_b}, and a constant polynomial reconstruction with no wet stencils is illustrated in Figure \ref{fig:Stencils_c}. 

\subsection{Discrete mass conservation} \label{sec:Discrete_Mass}
FV schemes are inherently mass-conserving, and this is a property we wish to maintain with the subgrid. Here, we ensure mass conservation by determining the free surface elevation such that the volume on both meshes is equal. Equating the discrete volumes on the two meshes gives 
\begin{equation}
    |\tilde{T}_m|\bar{h}_m = \sum_{k=1}^{N_m^{sg}} |T_{\ell(m,k)}| h_{\ell(m,k)}, \label{eq:Volume_Equality}
\end{equation}
where $\bar{h}_m$ is the integral averaged water depth defined in \eqref{eq:Int_Avg}.
We will start of by consider the first-order accurate case with $\mathbf{w}_m^{\eta}(\vec{\mathbf{x}}) = \bar{\eta}_m$.
The subgrid water depth is defined using the nonlinear max function; hence, determining the free surface directly is non-trivial. However, we can make it linear by only considering the wet subgrid cells. By using the definition of the water depth \eqref{eq:hsg}, we get 
\begin{equation}
    \sum_{k=1}^{N_m^{sg}} |T_{\ell(m,k)}| \max(0, \bar{\eta}_m + d_{\ell(m,k)}) = \sum_{k\in\mathcal{W}_m} |T_{\ell(m,k)}|(\bar{\eta}_m + d_{\ell(m,k)}), 
\end{equation}
where $\mathcal{W}_m$ is the wet set containing only the wet subgrid elements in the $m$'th coarse element
\begin{equation}
    \mathcal{W}_m = \{ k  ~|~  \bar{\eta}_m + d_{\ell(m,k)} > 0, ~ k = 1,2,...,N^{sg}_m\}.
\end{equation}
By using the wet set, we can find the free surface elevation that preserves the volume
\begin{equation}
    \bar{\eta}_m = (|\tilde{T}_{m}| \bar{h}_{m} - \sum_{k\in\mathcal{W}_m} |T_{\ell(m,k)}|d_{\ell(m,k)})(\sum_{k\in\mathcal{W}_m} |T_{\ell(m,k)}|)^{-1}, \label{eq:FS_sg}
\end{equation}
hence, the problem becomes finding the wet set $\mathcal{W}_m$. In general, we have that dry subgrid cells have zero water depth and the free surface is 
\begin{equation}
    \eta_{\ell(m,k)} = d_{\ell(m,k)}, \quad k \notin \mathcal{W}_m.
\end{equation}

In reality, when determining the free surface elevation, we distinguish between the cells' wet/dry state. If an element is completely wet -- all subgrid elements are included in the wet set $\mathcal{W}_m = \{k=1,2,...,n^{sg}_m \}$ -- the bathymetry in \eqref{eq:FS_sg} becomes the coarse mesh integral average and we get that
\begin{equation}
    \bar{\eta}_m = \bar{h}_{m} - \bar{d}_{m}.
\end{equation}
This is a direct consequence of the integral-averaged variables and the volume equality, which keeps the formulation consistent. 

If an element is partially wet, we use a Newton-Raphson method, to update the free surface elevation and check if the wet set is correct. The Newton-Raphson method is based on \eqref{eq:Volume_Equality}
\begin{equation}
    f(\bar{\eta}^{p}_m) = |\tilde{T}_m|\bar{h}_m  - \sum_{k=1}^{n_m^{sg}} |T_{\ell(m,k)}| \max(0, \bar{\eta}^p_m + d_{\ell(m,k)}), 
\end{equation}
where $p$ is the iteration count and the updated free surface is given by
\begin{equation}
    \bar{\eta}^{p+1}_m = \bar{\eta}^{p}_m - \frac{f(\bar{\eta}^{l}_m)}{f^{\prime}(\bar{\eta}^{l}_m)}. 
\end{equation}

If the element is entirely dry $\mathcal{W}_m = \varnothing$, the free surface is equal the bathymetry in every subgrid cell.  

To maintain consistency with the integral-averaged variables, the integral-averaged bathymetry is recomputed based on the wet subgrid cells. The reevaluated bathymetry can be computed using the wet set or the newly computed free surface elevation
\begin{equation}
    \bar{d}_m^{wet} = \bar{h}_m - \bar{\eta}_m. 
\end{equation}
This is, of course, only necessary for cells that are either partially wet or completely dry, as wet cells fulfill this requirement already. For simplicity, the '$wet$' superscript is dropped, and it is assumed for the remainder of the paper.

\subsubsection{Mass conservation for the second-order reconstruction}
For the second-order case, we begin by constructing the free surface elevation by assuming it is constant, as just described by the discrete mass algorithm. The polynomials are adjusted to ensure the volume equality is fulfilled. The subgrid water depth is given as
\begin{equation}
    h_{\ell(m,k)} = \max(0, \bar{\eta}_m + \nabla\eta(\vec{\mathbf{x}}_{\ell(m,k)} - \vec{\mathbf{x}}_m) + d_{\ell(m,k)}).
\end{equation}
In order to fulfill the volume criteria without recomputing the polynomial, we simply freeze the gradients and shift the polynomial by changing $\bar{\eta}_m$ iteratively. For this, we use the same Newton-Raphson method as for the constant case.

\subsubsection{Note on water distribution algorithm}
The current water distribution algorithm exhibits two notable shortcomings: non-physical transport of water across unresolved barriers and unfavorable water redistribution on steep slopes. The presented free surface algorithm can been seen as a cell-local bathtub approach, as it redistributes the water volume from the deepest subelement until the water volumes are balanced. If a physical barrier is placed in the middle of an element, it is ignored by the algorithm, and water may be transported across the barrier.  

A similar issue arises on steep slopes, where water is distributed from the deepest subelement, resulting in most of the water being concentrated at one end of the element and the flow being modeled as a series of waterfalls. This issue becomes more pronounced as the number of subelements within a coarse cell increases. However, it should be noted that this is only problematic when the flow is parallel to the slope. If the flow is perpendicular to the slope, the approach may be advantageous because it helps preserve the flow direction, for example in a river channel constrained by dikes.
Adjustments to mitigate these issues have been proposed in
\cite{Begmohammadi2021} and \cite{STELLING2022}.

\subsection{Face value reconstruction} \label{sec:ValRec}
The values at the faces are important for the Riemann problem and for getting a well-balanced scheme. Here we use some of the key ideas from \cite{Chen2017}, \cite{Audusse2004}, and \cite{Shamkhalchian2021}.

\begin{figure}[htb]
    \centering
    \vspace{3mm}
    \subfloat[Case 1.1]{
  \includegraphics[width=65mm]{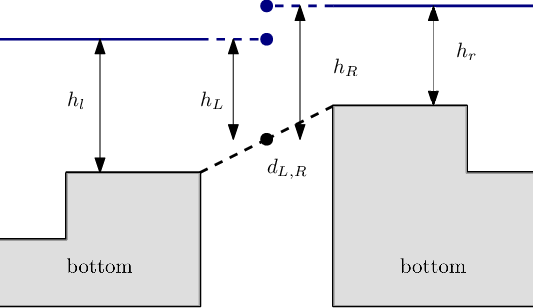}
}   \vspace{3mm}
    \hspace{5mm}
    \subfloat[Case 1.2.1]{
  \includegraphics[width=65mm]{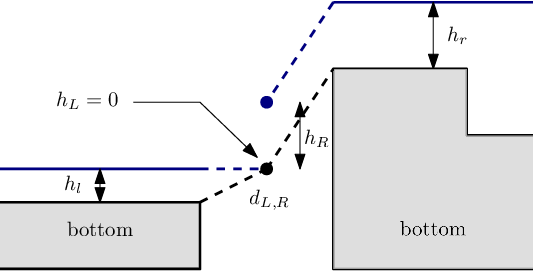} \label{fig:FV_C121}
} 
    \hspace{5mm}
    \subfloat[Case 1.2.2]{
  \includegraphics[width=65mm]{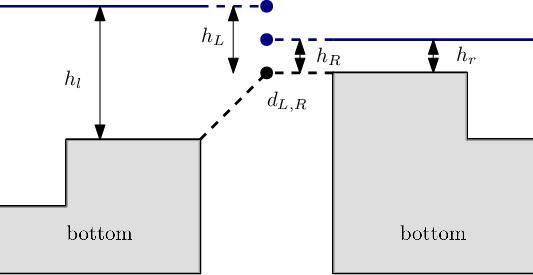} \label{fig:FV_C122}
}   \vspace{4mm}
    \hspace{5mm}
    \subfloat[Case 2.1]{
  \includegraphics[width=65mm]{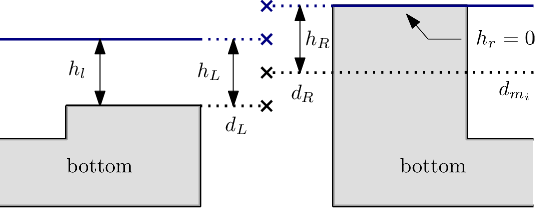}
} 
    \hspace{5mm}
    \subfloat[Case 2.2]{
  \includegraphics[width=65mm]{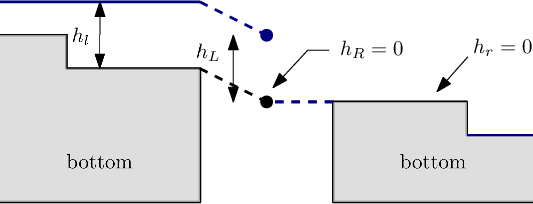}
}
    \hspace{5mm}
    \subfloat[Case 3]{
  \includegraphics[width=65mm]{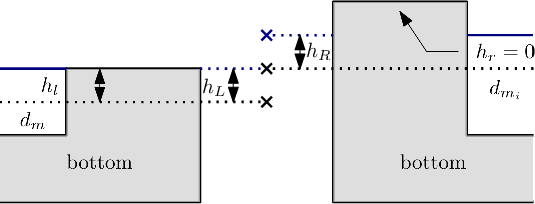} \label{fig:FV_C3}
}
    \caption{Illustrations of 6 different face value reconstruction cases and how it is treated. The dashed reconstructions with circles are for the Riemann problem, and the dotted with crosses are for the wall flux.} \label{fig:FaceVal_recon}
\end{figure}

When reconstructing the values at the faces, we distinguish between three different cases, each with its own individual subcases. The cases are based on the values from the subgrid cells sharing a face. The three cases are that: 1) both subgrid cells are wet, 2) one is wet and the other is dry, and 3) both are dry. Note that the last case is possible due to the partially wet cells. All cases are illustrated in Figure \ref{fig:FaceVal_recon}. \\
In order to simplify the notation we use a left and right notation where $l$ is the extrapolated value from left and $r$ from the right to the same face. As previously mentioned, $l$ is always cell $T_{\ell(m,k)}$ and $r$ is the neighbor across face $j$.   
To clarify the values at the faces are 
\begin{equation}
    d_l = d_{\ell(m,k)}, \quad \eta_l = \mathbf{w}_m^{\eta}(\vec{\mathbf{x}}_{\ell(m,k),j}), \quad \text{and} \quad h_l = \max(0,\mathbf{w}_m^{\eta}(\vec{\mathbf{x}}_{\ell(m,k),j})+d_{\ell(m,k)}),
\end{equation}
where $\vec{\mathbf{x}}_{\ell(m,k),j}$ are the coordinates of the subface center. In a similar fashion the right state can be determined.
After the reconstruction, we refer to it by $L$ and $R$, which are the values used for the numerical flux. 

Case 1, where both sides are wet. Here we distinguish between whether the water height on either side is less than the difference in bathymetry. The water height condition is given by 
\begin{equation}
    h_l < |d_l - d_r| \quad \text{or} \quad h_r < |d_l - d_r|. 
\end{equation}
If the condition is not met, we use
\begin{subequations}
\begin{align}
    &d_{L,R} = \frac{1}{2}(d_l + d_r), \\ 
    &h_L = \max(\eta_l + d_L ,0), \\ 
    &h_R = \max(\eta_r + d_R ,0),
\end{align}    
\end{subequations}
where $d_{L,R} = d_L = d_R$ is the notation for the bathymetry reconstruction being the same for both sides of the face. This is considered Case 1.1.

If the condition is met we use the approach from \cite{Chen2017} given as
\begin{subequations}
    \begin{align}
    &d_{L,R} = - \min(-\min(d_l,d_r),\min(\eta_l,\eta_r)), \\ 
    &h_L = \min(\eta_l + d_L ,h_l), \\ 
    &h_R = \min(\eta_r + d_R ,h_r).
    \end{align}    
\end{subequations}
This is considered Case 1.2. Depending on the initial extrapolated values the reconstruction can look quite differently, see e.g. Figure \ref{fig:FV_C121} and \ref{fig:FV_C122}. 

Case 2, with one wet and one dry cell, involves either the left cell being wet and the right cell being dry, or vice versa; both cases are treated identically. Here, we focus on the case where the left cell is wet and the right cell is dry. We have two different subcases, one if the bathymetry on the right is higher than the free surface on the left ($\eta_l < -d_r$). In this case, we do not expect water to be transported between the cells, and we treat it as a wall, with 
\begin{subequations}
\begin{align}
    &d_L = d_l, \\
    &h_L = \eta_l + d_l.
\end{align}
\end{subequations}
For the right side, we use the coarse mesh variables
\begin{subequations}
\begin{align}
    &d_R = \bar{d}_{m_i}, \\ &h_R = \bar{\eta}_{m_i} + \bar{d}_{m_i}. 
\end{align}
\end{subequations}
This is considered Case 2.1.

If the free surface is above the neighbor's bathymetry, we expect transport from left to right and proceed to Case 1).
This is considered Case 2.2.

Case 3, with both sides being dry, we set the bathymetry and the water depth to be the coarse mesh values
\begin{subequations}
\begin{align}
    &d_L = \bar{d}_{m}, \\
    &d_R = \bar{d}_{m_i}, \\ 
    &h_L = \bar{h}_{m} = \bar{\eta}_{m} + \bar{d}_{m}, \\ 
    &h_R = \bar{h}_{m_i} = \bar{\eta}_{m_i} + \bar{d}_{m_i}.
\end{align}
\end{subequations}
This is treated as a wall -- like Case 2.1 -- as it represents no water transport across the face. 

\subsection{Numerical flux treatment} \label{sec:Flux}
The numerical flux is treated differently depending on the face value reconstruction case, with either a solution to the Riemann problem or a wall condition.

To solve the Riemann problem, we use the approximate HLLC solver by E. Toro, including the dry solution. The HLLC solver is documented extensively in the literature and will not be outlined here; for a summary, see \cite{Fraccarollo1995}, \cite{Toro2001}, or \cite{Toro2019}. The HLLC solver is used for Case 1 and Case 2.1. 

The Riemann solver computes a unified flux across the face. For the gravity source term, we extract a unified water depth and consequently the free surface elevation, as 
\begin{equation}
    h^*_{\ell(m,k),j} = h^{\text{HLLC}}, \quad \eta^*_{\ell(m,k),j} = h^{\text{HLLC}} - d_{L,R}, \quad \text{and} \quad d^*_{\ell(m,k),j} = d_{L,R},
\end{equation}
where $*$ denotes the values are extracted from the flux calculation and with details as to how $h^{\text{HLLC}}$ is computed is provided in Appendix \ref{sec:appendix}.
Likewise, we have the face values for the wall condition, which are unchanged from the face value reconstruction cases. They are 
\begin{equation}
    h^*_{\ell(m,k),j} = h_{L}, \quad \eta^*_{\ell(m,k),j} = \eta_{L},  \quad \text{and} \quad  d^*_{\ell(m,k),j} = d_{L}.
\end{equation}

The wall condition is a zero-velocity impermeability condition, $\mathbf{n}\cdot\mathbf{u} = 0$, which reflects the incoming momentum. When inserting this into the flux \eqref{eq:flux}, we get the momentum flux given purely by the gravity part
\begin{equation}
    \mathcal{F}^{\text{wall}}(\bar{\mathbf{U}}_L,\bar{\mathbf{U}}_R,\mathbf{n}) = \begin{bmatrix}
        0 \\
        \frac{1}{2}g((\eta^*_{\ell(m,k),j})^2 + 2\eta^*_{\ell(m,k),j} d^*_{\ell(m,k),j})n^x_{m,i} \\
        \frac{1}{2}g((\eta^*_{\ell(m,k),j})^2 + 2\eta^*_{\ell(m,k),j} d^*_{\ell(m,k),j})n^y_{m,i}
    \end{bmatrix},
\end{equation}
where the wall values are reconstructed from cell $T_{\ell(m,k)}$. It is important to note that in these cases, the flux from the left and the right side is not necessarily the same. 
This is applied in reconstruction Case 2.1 and Case 3, which are cases with at least one dry side.

\subsection{Gravity source term} \label{sec:Gravity}
We present a new discretization of the gravity source term constructed to balance the gravity flux and source term on unstructured meshes. We start off by stating the source term discretization on the coarse and subgrid mesh, where a thorough derivation is given in Appendix \ref{sec:appendix}. After, we show that it is indeed a valid approximation and that it is well-balanced.

The gravity source term is transformed to a line integral over the element faces and are given on the coarse mesh as
\begin{equation}
    \mathcal{S}_{g}(\bar{\mathbf{U}}_m) = \frac{1}{|\tilde{T}_m|}\frac{1}{2}g\sum_{i=1}^3(\eta^*_{m,i} + \bar{\eta}_m)(d^*_{m,i} - \bar{d}_m) |\tilde{E}_{m,i}| \mathbf{n}_{m,i}, \label{eq:gravity_C}
\end{equation}
The gravity source term is extended to the subgrid by summing over the subgrid edges 
\begin{equation}
    \mathcal{S}_{g}(\bar{\mathbf{U}}_m) = \frac{1}{|\tilde{T}_m|}\frac{1}{2}g \sum_{i=1}^3 \left( \sum_{(k,j)\in\mathcal{E}_{m,i}} (\eta^*_{\ell(m,k),j} + \bar{\eta}_m)(d^*_{\ell(m,k),j} - \bar{d}_m)|E_{\ell(m,k),j}|\mathbf{n}_{m,i}\right). \label{eq:gravity_sg}
\end{equation}

\subsubsection{Discretization properties}
To show that \eqref{eq:gravity_C} is valid, we expand the terms and use the inverse discrete divergence theorem
\begin{equation}
    \int_{\tilde{T}_m} \mathbf{S}_g(\mathbf{U}) ~\text{d}\vec{\mathbf{x}} = \frac{1}{2}g\int_{\tilde{T}_m} \nabla\eta d - d\nabla\eta + \eta \nabla d~\text{d}\vec{\mathbf{x}},
\end{equation}
where the line integral is transformed back into a surface integral. 
Now, applying the product rule, most of the derivatives cancel each other, and we end up with
\begin{equation}
    \int_{\tilde{T}_m} \mathbf{S}_g(\mathbf{U}) ~\text{d}\vec{\mathbf{x}} = \int_{\tilde{T}_m} g\eta \nabla d ~\text{d}\vec{\mathbf{x}},
\end{equation}
which is identical to the original formulation in \eqref{eq:source}. 

Another property of this discretization is that we can show the scheme is well-balanced in the sense of \cite{LEVEQUE1998} for the initial state 
\begin{equation}
    \eta(\vec{\mathbf{x}},0) = \eta_0,   \quad \vec{\mathbf{u}}(\vec{\mathbf{x}},0) = 0, \quad \text{where} \quad d(\vec{\mathbf{x}}) \neq const,
\end{equation}
and $\eta_0 = const$. The discharge vanishes for $\vec{\mathbf{u}}(\vec{\mathbf{x}},0) = 0$, hence only the balancing of the gravity terms remain $\mathcal{S}_{g}(\bar{\mathbf{U}}_m) - \mathcal{F}_{g}(\bar{\mathbf{U}}_m) = 0$. Both the flux term and the source term are rewritten to be dependent on water depth and bathymetry alone
\begin{subequations}
\begin{align}
        |\tilde{T}_m|\mathcal{F}_{g}(\bar{\mathbf{U}}_m) = \frac{1}{2}g\sum_{i=1}^3[(h^*_{m,i} + \bar{h}_m)(h^*_{m,i} - \bar{h}_m) - (d^*_{m,i} + \bar{d}_m)(d^*_{m,i} - \bar{d}_m)]|\tilde{E}_{m,i}|\mathbf{n}_{m,i}, \\
        |\tilde{T}_m|\mathcal{S}_{g}(\bar{\mathbf{U}}_m) = \frac{1}{2}g\sum_{i=1}^3[(h^*_{m,i} + \bar{h}_m)(d^*_{m,i} - \bar{d}_m) -  (d^*_{m,i} + \bar{d}_m)(d^*_{m,i} - \bar{d}_m)]|\tilde{E}_{m,i}|\mathbf{n}_{m,i}.
\end{align}
\end{subequations}
When subtracting the terms, the bathymetry part cancels, and by rearranging the terms, we get that 
\begin{equation}
    |\tilde{T}_m|(\mathcal{S}_{g}(\bar{\mathbf{U}}_m) - \mathcal{F}_{g}(\bar{\mathbf{U}}_m)) = \frac{1}{2}g\sum_{i=1}^3(h^*_{m,i} + \bar{h}_m)(\bar{\eta}_m-\eta^*_{m,i})|\tilde{E}_{m,i}|\mathbf{n}_{m,i},
\end{equation}
which is balanced when $\bar{\eta_m} = \eta^*_{m,i}$. From the initial condition, we have that the free surface is constant, and therefore, the scheme is well-balanced.
Now, we consider the subgrid gravity source term \eqref{eq:gravity_sg} and the two relevant states; completely wet and partially wet. It follows directly from the coarse mesh that the expression is still well-balanced for completely wet cells. For partially wet cells, we use the reconstruction from Section \ref{sec:ValRec}, where we have that $\eta^*_{m,i} = \bar{\eta}_m$ for Case 2.1 and Case 3. With this, the scheme is well-balanced, also with the subgrid. 

The gravity source term \eqref{eq:gravity_sg} also reduces computational complexity, as the surface integral is transformed into a line integral. For the subgrid method, this means that it runs over $3\cdot n^{sg}_m$ faces instead of $(n^{sg}_m)^2$ subgrid elements.  

\subsection{Friction source term}\label{sec:Fric_Coef}
The friction term has a nonlinear dependency on both the bathymetry and roughness, and can lead to large errors if the resolution is insufficient \cite{Volp2013}. In the article, they showed that by incorporating subgrid, the errors from linearization are reduced. 

Unlike other explicit subgrid methods, we do not employ the ideas of \cite{Volp2013}, as we found that the assumption of matching bathymetry and free surface gradient is not general enough.    
Here, we use a more straightforward averaging approach that stays closer to the original formulation. Recall that the friction term is given by  
\begin{equation}
    \mathcal{S}_{f}(\bar{\mathbf{U}}_m) = \frac{1}{|\tilde{T}_m|}\int_{\tilde{T}_m} \vec{\mathbf{u}} |\vec{\mathbf{u}}|c_{fb} ~\text{d} \vec{\mathbf{x}}.
\end{equation}
The friction term is difficult to represent in partially wet cells as it is undefined for dry subgrid cells. We found that averaging over the wet subgrid cells and using a constant velocity -- even for second-order -- works well. The friction term becomes
\begin{equation}
    \mathcal{S}_{f}(\bar{\mathbf{U}}_m) = \vec{\mathbf{u}}_m |\vec{\mathbf{u}}_m| c_{fb,m}^{sg},
\end{equation}
where the discrete friction coefficient is based on wet subgrid cells only and considers subgrid water depth and Manning's $M$
\begin{equation}
    c_{fb,m}^{sg} = \left(\sum_{k\in\mathcal{W}_m} |T_{\ell(m,k)}|\right)^{-1} \sum_{k\in\mathcal{W}_m}\frac{g}{M^2_{\ell(m,k)}h^{1/3}_{\ell(m,k)}}|T_{\ell(m,k)}|, \label{eq:cfbsg}
\end{equation}

\subsection{Temporal integration}

For temporal integration, we use a semi-explicit TVD Runge-Kutta scheme based on an explicit scheme (see e.g. \cite{gottlieb1998total})
\begin{subequations}
    \begin{align}
        \bar{\mathbf{U}}_m^{(1)} &= \bar{\mathbf{U}}_m^{n} + \Delta t (\hat{\mathcal{F}}(\bar{\mathbf{U}}^n_m)+\mathcal{S}_g(\bar{\mathbf{U}}^n_m) + \mathcal{S}_f(\bar{\mathbf{U}}^{(1)}_m)),\label{eq:RK2a}\\
        \bar{\mathbf{U}}^{n+1}_m &= \frac{1}{2}\bar{\mathbf{U}}^n_m + \frac{1}{2}\left[\bar{\mathbf{U}}_m^{(1)} + \Delta t(\hat{\mathcal{F}}(\bar{\mathbf{U}}^{(1)}_m)+\mathcal{S}_g(\bar{\mathbf{U}}^{(1)}_m) + \mathcal{S}_f(\bar{\mathbf{U}}^{n+1}_m)) \right],
    \end{align} \label{eq:RK2}
\end{subequations}
where superscript $(1)$ is the stage and $n$ is the current time step number. The scheme is second-order accurate and for the first-order version we use the initial Euler-step only \eqref{eq:RK2a}. The velocities in the first stage and the time step $(n+1)$ are corrected using an implicit friction correction. For the time step, only the expression in the square brackets is corrected. By doing the correction on part of the time step, we are able to preserve second-order for an implicit correction that is otherwise first-order. 

The implicit method is outlined in \cite{Xia2017}, which proposed a Newton-Raphson method for correcting the velocities. For better conditioning of the system, they suggest using the water depth at time $n+1$ to determine the velocities $\vec{\mathbf{u}}_{*,m} = \vec{h\mathbf{u}}_m^n/\bar{h}_m^{n+1}$. The scheme is given by
\begin{equation}
    \vec{\mathbf{u}}_{*,m}^{p+1} = \vec{\mathbf{u}}_{*,m}^{p} + (\mathbf{I} - \Delta t \mathcal{J}_S)^{-1} [\Delta t \mathcal{S}_f(\
    \vec{\mathbf{u}}_{*,m}^{p}) + 
    \frac{\Delta t}{\bar{h}_m^{n+1}}(\mathcal{F}(\bar{\mathbf{U}}_m^n) + \mathcal{S}_g(\bar{\mathbf{U}}_m^n)) - (\vec{\mathbf{u}}_{*,m}^{p} - \vec{\mathbf{u}}_{*,m})],
\end{equation}
where $p$ is the number of the Newton iteration, $\mathbf{I}$ is the identity matrix, and $\mathcal{J}_S = \partial \mathcal{S}_f / \partial \vec{\mathbf{u}}$ is the Jacobian of the friction source term. 
The iterations are carried out on the coarse mesh, even though the friction term is evaluated on the subgrid. This is because the friction coefficient $c_{fb}^{sg}$ in \eqref{eq:cfbsg} depends only on the already computed water depth and is therefore invariant during the iterations. Consequently, the correction of the friction term depends solely on the velocities defined on the coarse mesh.

\subsubsection{Stability condition} \label{sec:CFL}
The conserved variables are given on the coarse mesh; hence, so is the domain of dependence. This means that the stability condition is governed by the coarse mesh. Here, we use a classic CFL condition, see e.g. \cite{BOSCHERI2023,Yoon2004} 
\begin{equation}
    \Delta t \leq CFL \min_{m=1,2,...,N_c}\left(\frac{\Delta s_m}{\lambda_m}\right), \label{eq:Delta_t}
\end{equation}
where $\Delta s_m = \sqrt{|\tilde{T}_m|}$ is the characteristic length and $\lambda_m$ is the eigenvalue corresponding to the normal wave speed. The CFL conditions should be $<\frac{1}{2}$ to ensure stability on 2D unstructured meshes. The eigenspeed is given by the celerity and the velocity  
\begin{equation}
    \lambda_m = |\vec{\mathbf{u}}_m| + \sqrt{gh_{\ell(m,k)}^{max}},
\end{equation}
where $\displaystyle h_{\ell(m,k)}^{max}= \max_{k=1,2,..,N^{sg}_m}(h_{\ell(m,k)})$ is the largest subgrid water depth. We use the subgrid water depth because $\bar{h}_m \leq h_{\ell(m,k)}^{max}$. This implies that the time step will be smaller or equal, making it a conservative choice. Similar to other well-balanced SWE schemes with friction, the stability criterion is based solely on the hyperbolic part of the equations \cite{Xia2017,Chertock2015}. The friction term has a damping effect on the solution, making this choice mildly conservative.

\subsection{Model complexity} \label{sec:Model_scale}
The equal subgrid division used here results in a total number of $N_m^{sg}=(n_m^{sg})^2$ subgrid elements for each coarse element. The total number of subgrid faces overlapping the coarse faces is $3\cdot n_m^{sg}$. Since flux and gravity calculations are performed on the faces rather than on each subgrid element, the subgrid method quickly becomes attractive, as it requires only $N_c\cdot3\cdot n^{sg}_m$ computations. In comparison, a traditional method that resolves every subgrid cell explicitly requires $3 \cdot N_c\cdot N_m^{sg}$ computations. This corresponds to a reduction by a factor of $n_m^{sg}$ per time stage. 
The time integration is independent of the subgrid discretization and therefore introduces no additional computations. Furthermore, the CFL condition is primarily determined by the coarse mesh, meaning that the number of time steps is unaffected by the number of subgrid elements.

Unfortunately, some operations depend on all subgrid elements within a coarse element. In particular, the free-surface algorithm and the friction coefficient calculations both require information from every subgrid element. However, we have found that the Riemann problem remains the computational bottleneck, at least for the subgrid divisions considered here and in the current implementation, which has not yet been fully optimized. 

To make the model scalable with respect to the number of subface divisions $n_m^{sg}$, lookup tables can be used, as suggested in \cite{GUINOT2018,Sanders2019}. It should be noted that this approach is only applicable to the first-order free surface algorithm. For the friction coefficient calculation, this would turn it into a piecewise-linear approximation.

%% file: Texts/Results_verification.tex
\section{Model Verification and Validation} \label{sec:Model_Verfication}
This section presents results obtained with the first- and second-order schemes. The goal is to present the correctness of the subgrid method, which is why we mainly use test cases with exact analytical solutions. Many of the cases are from the collection of solutions given in \cite{Delestre2013}, with cases containing friction and moving flood and dry boundaries. For some of the cases we also use no subgrid -- referred to as \textit{Classic} -- for better comparison, where the second-order no subgrid scheme is true second-order with bathymetry gradients. 

All coarse meshes are 2D unstructured triangular meshes generated with Gmsh \cite{Geuzaine2009}. Boundary conditions are applied at the subgrid level following the approach of \cite{Sleigh1998}. In all test cases, the stability condition from \eqref{eq:Delta_t} is used to determine $\Delta t$ with $CFL = [0.4;0.5]$, unless otherwise stated. The dry tolerance for the cut-off velocities are $\epsilon_{dry} = 10^{-4}$ m. Errors are measured in the $L_{\infty}$-norm and the $L_2$-norm 
% \begin{subequations}
% \begin{align}    
%     L_{\infty}(h) = \max_{\Omega}(|\mathbf{w}^h_m(\vec{\mathbf{x}}) - h_e(\vec{\mathbf{x}})|), \\ L_{2}(h) = \sqrt{\int_{\Omega}(\mathbf{w}^h_m(\vec{\mathbf{x}}) - h_e(\vec{\mathbf{x}}))^2~\text{d}\vec{\mathbf{x}}},
% \end{align}
% \end{subequations}
\begin{subequations}
\begin{align}    
    L_{\infty}(h) = \max_{\Omega}(|h(\vec{\mathbf{x}}) - h_e(\vec{\mathbf{x}})|), \\ 
    L_{2}(h) = \sqrt{\int_{\Omega}(h(\vec{\mathbf{x}}) - h_e(\vec{\mathbf{x}}))^2~\text{d}\vec{\mathbf{x}}},
\end{align}
\end{subequations}
where the water depth $h$ is used as an example, and where the subscript $e$ denotes the exact solution. The norms are computed on the coarse mesh in the cell centers for easy comparison between subgrid and no subgrid.   

\subsection{Well-balanced test}
To verify that the scheme is well-balanced and fulfills the C-Property, we follow the benchmark proposed in \cite{LEVEQUE1998}, also given in \cite{Delestre2013}. The test consists of two examples without friction; one to test the scheme's ability to numerically preserve a stationary equilibrium, and two to test that spurious oscillations do not occur in the presence of a small perturbation. 

\begin{figure}[htb]
    \centering
    \includegraphics[width=0.55\linewidth]{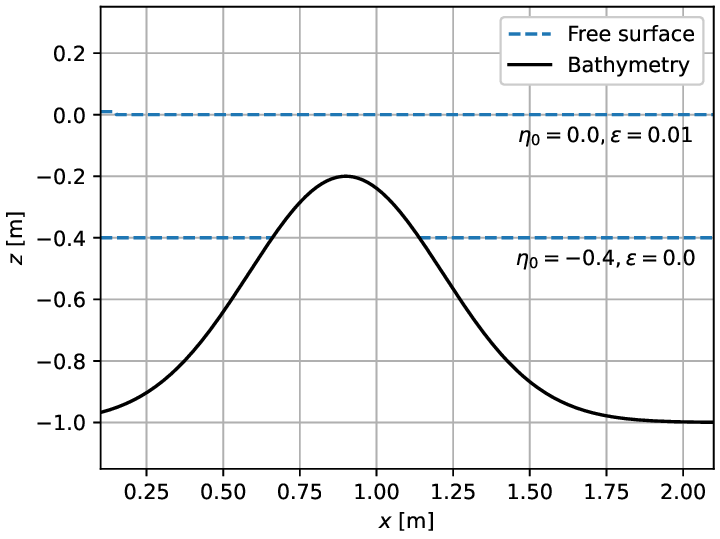}
    \caption{Cross-sectional view at $y = 0.5$ m of the bathymetry and the initial free surface elevations for the two test cases.}
    \label{fig:WB_IC}
\end{figure}

The computational domain is $\Omega = [0.1;2.1] \times [0;1]$ m\textsuperscript{2} and is equipped with wall boundary conditions ($\mathbf{n}\cdot\vec{\mathbf{u}} = 0$) in the $x$-direction and periodic boundary conditions in the $y$-direction. The domain is discretized with a coarse mesh having 11,628 coarse elements with a side length of $|\tilde{E}_m| \approx 1/50$ m. For the subgrid mesh two different ones are considered, namely a subgrid division of $n^{sg}_m = 2$ and $n^{sg}_m = 5$.  
The bathymetry is given by
\begin{equation}
    d(\vec{\mathbf{x}}) = 1 - 0.8e^{-5(x-0.9)^2-50(y-0.5)^2},
\end{equation}
and the initial state is
\begin{equation}
    \eta(\vec{\mathbf{x}},0) = \begin{cases}
    \eta_0 + \varepsilon \quad &\text{if }~ x \leq 0.15 \\
    \eta_0 & \text{elsewhere}
    \end{cases} \quad \text{and} \quad \vec{\mathbf{u}}(\vec{\mathbf{x}},0) = 0,
\end{equation}
where $\eta_0$ is the constant free surface elevation level and $\varepsilon$ is the amplitude of the perturbation, see Figure \ref{fig:WB_IC}. 

\subsubsection{Lake at rest with emerged bump}
In order to show that water at rest is preserved, the initial state is set to  $\eta_0=-0.4$ m and $\varepsilon = 0.0$ m such that the top of the bathymetry is above the free surface elevation. This creates a flood and dry boundary around the emerged bump, demonstrating that the scheme can preserve zero velocity in the presence of partially wet cells. The errors are recorded at the final time $t_f = 0.1$ s after 100 fixed time steps of $\Delta t = 0.001$ s. The errors are reported in Table \ref{tab:Lake_at_Rest}, where all results are more or less machine precision.  

\begin{table}[htb]
    \centering
    \caption{Well-balanced test with different polynomial order and number of subgrid cells. Errors are measured in $L_\infty$- and $L_2$- norms for $h$ and $hu$ at the final time $t_f = 0.1$ s.}
    \begin{tabular}{cccccc}
         & $n^{sg}_m$ & $L_\infty(h)$ & $L_2(h)$ & $L_\infty(hu)$ & $L_2(hu)$ \\ \hline
     \multirow{2}{*}{$P_n = 0$} 
         & 2 & 1.875E-13 & 1.751E-16 & 1.018E-13 & 1.865E-15 \\ 
         & 5 & 3.369E-13 & 1.751E-16 & 1.153E-13 & 1.773E-15 \\ \hline
     \multirow{2}{*}{$P_n = 1$}    
        & 2 & 3.174E-13 & 7.264E-16 & 4.520E-13 & 2.417E-15\\
        & 5 & 2.199E-13 & 7.228E-16 & 2.993E-13 & 2.363E-15\\ \hline
    \end{tabular}    
    \label{tab:Lake_at_Rest}
\end{table}

\subsubsection{Small perturbation}
\begin{figure}[htb]
    \centering
    % \vspace{3mm}
    \subfloat[t = 0.12 s]{
    \includegraphics[width=74mm]{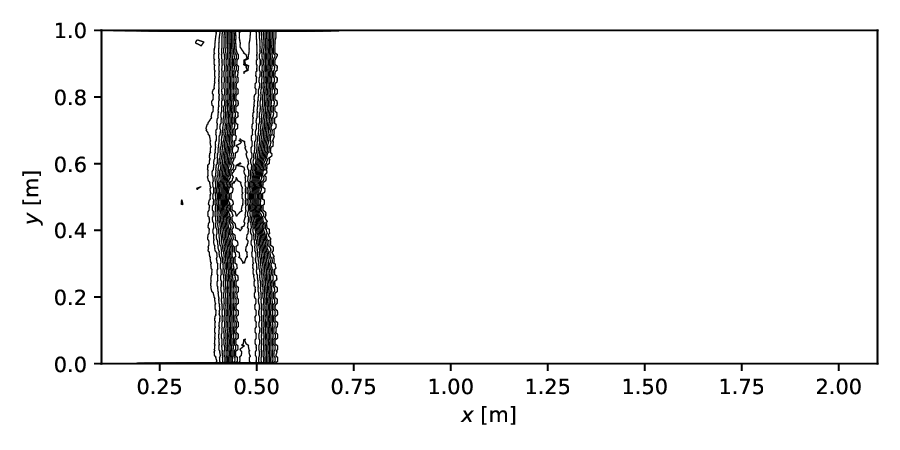}
}   
    \subfloat[t = 0.24 s]{
    \includegraphics[width=74mm]{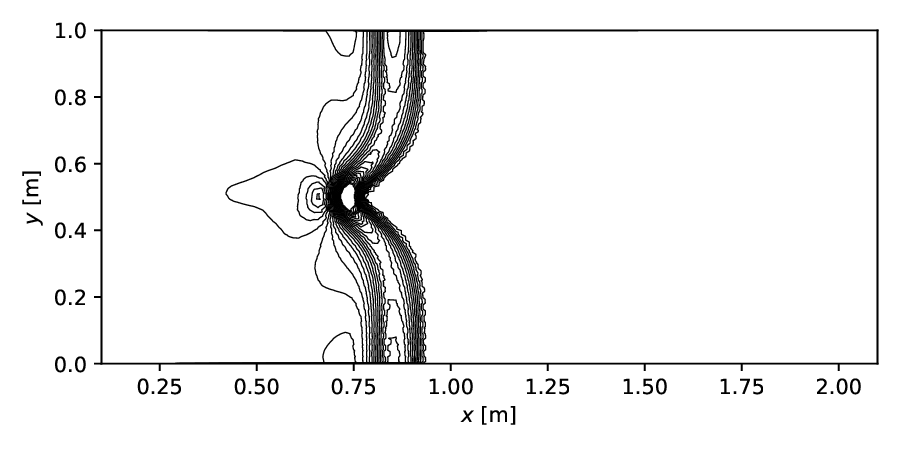} 
}

    \subfloat[t = 0.36 s]{
    \includegraphics[width=74mm]{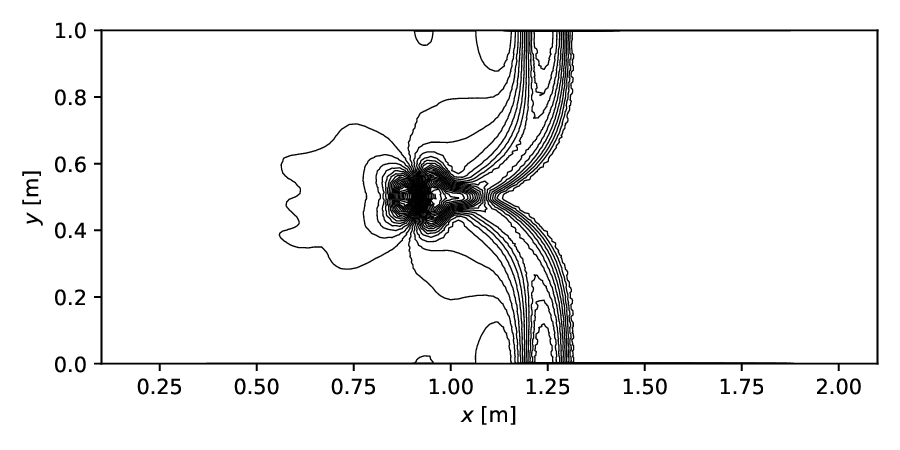}   
}
    \subfloat[t = 0.48 s]{
    \includegraphics[width=74mm]{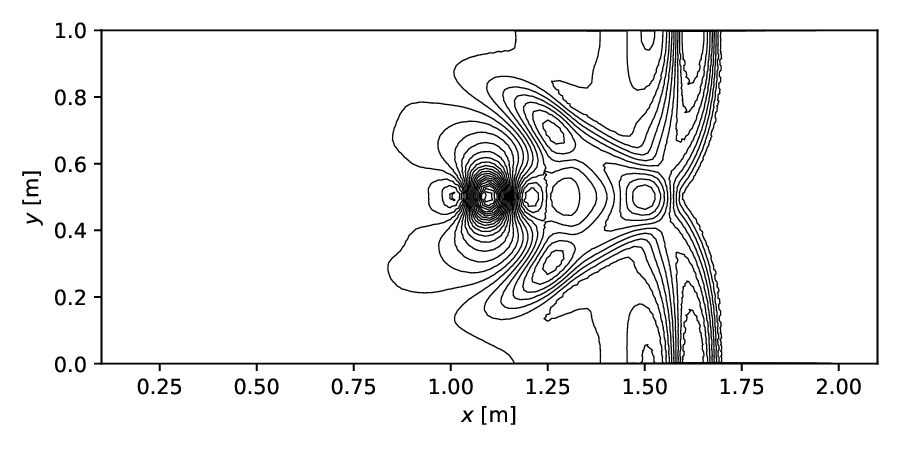}
} 

    \caption{Results of small perturbation of the free surface for the second-order scheme. This is illustrated with 40 equidistant contour lines in the interval $\eta=[0.993; 1.008]$ m. } \label{fig:Small_Pertubation}
\end{figure}

A small perturbation is added to the free surface elevation, with parameters $\eta_0 = 0.0$ m, $\varepsilon = 0.01$ m, and a final time of $t_f = 0.48$ s with variable time stepping. Here, the bump is totally immersed in water, making the domain entirely wet. Results of the evolution of the free surface elevation are given in Figure \ref{fig:Small_Pertubation} at different times. The scheme generates no spurious oscillations when the small perturbation travels over the bump. 
The results are in agreement with other results reported in the literature \cite{XING2006,CANESTRELLI2010,RICCHIUTO2015} and, of course, the original work \cite{LEVEQUE1998}. 

\subsection{Steady state solution over variable bed}
The order of accuracy with friction is verified with a MacDonald type test \cite{MacDonald1997}, which is a 1D steady state solution over a variable bed. The solution assumes constant discharge and that the solution attains a steady state, as a balance between gravity and friction. The solution to the problem is manufactured, and the bathymetry is found as a solution to the momentum equation in the $x$-direction (\ref{eq:gov_eq}-\ref{eq:source}) with a constant discharge. The bathymetry is found with a very fine discretized Runge-Kutta 4 method. 
The solutions were originally proposed by \cite{MacDonald1997}; here, we follow an example from \cite{Delestre2013}.

The channel is $L = 5000$ m long and since the solution is purely in the $x$-direction we can pick the width freely. The computational domain is $\Omega = [0;5000] \times [0;W_i]$, m\textsuperscript{2} with $W_i\in\{100,200,300,400\}$ m depending on the coarse element size. The friction is based on a constant Manning number of $M=\frac{1}{0.03}$ s/m\textsuperscript{1/3} in the entire domain. The domain is equipped with a discharge boundary condition in the west end of $\mathbf{n}\cdot h\vec{\mathbf{u}} = 2$ m\textsuperscript{2}/s, a constant water depth boundary condition in the east end of $h = h_{e}(L)$, and periodic boundary conditions in the $y$-direction. 

\begin{figure}[htb]
    \centering
    \includegraphics[width=0.65\linewidth]{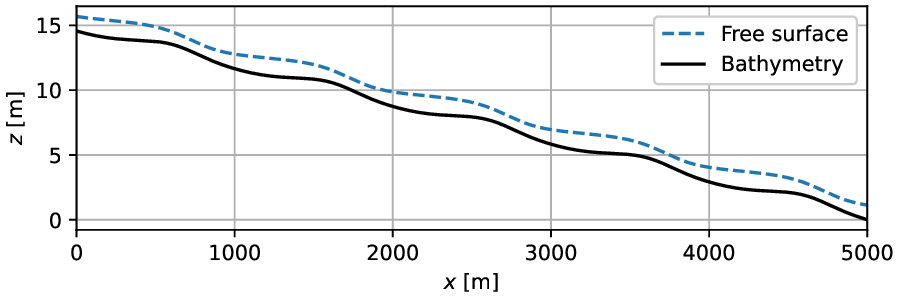}
    \caption{Steady state solution of the MacDonald type example, showing the free surface elevation and the bathymetry.}
    \label{fig:BM_McD5000}
\end{figure}

The manufactured water depth is
\begin{equation}
    h_{e}(x) = \frac{9}{8} + \frac{1}{4}\sin\left(\frac{\pi x}{500}\right).
\end{equation}
See Figure \ref{fig:BM_McD5000} for an illustration of the free surface and the bathymetry at steady state. 
The problem is initialized with a dry bed except for a little lake downstream, and zero velocity everywhere
\begin{equation}
    h(\vec{\mathbf{x}},0) = \max(h_{e}(L) - d(L) + d(x),0), \quad \text{and} \quad \vec{\mathbf{u}}(\vec{\mathbf{x}},0) = 0.
\end{equation}
The simulations are run until $t_f = 100$ minutes where the simulation has reached a steady state.  

In order to get an idea of the influence of the subgrid, a convergence study is made with a constant subgrid face division of $n_m^{sg}=$ 2, 4, and 8, and a decreasing coarse cell size of $|\tilde{E}_m| \approx$ 100, 75, 50, and 25 m. The errors are measured in a normalized $L_2$-norm ($L_2(\cdot)/|\Omega|$), as the width of the domain is adjusted with the coarse mesh. The $L_2$-error with respect to water depth and discharge, and the order of convergence can be seen in Table \ref{tab:MacDonald_Conv}.
% The $L_2$-error with respect to water depth and discharge, and the order of convergence can be seen in Figure \ref{fig:convL2h} and \ref{fig:convL2hu}, respectively.
% \begin{figure}[htb]
%     \centering
%     % \vspace{3mm}
%     \subfloat[Water depth $h$]{    \includegraphics[width=0.48\linewidth]{Figures/Conv_L2hV1.eps} \label{fig:convL2h}
% }
%     \subfloat[Discharge $hu$]{
%   \includegraphics[width=0.48\linewidth]{Figures/Conv_L2huV1.eps} \label{fig:convL2hu}
% }

%     \caption{Convergence study of a Macdonald type test with no subgrid (\textit{Classic}) and different numbers of subgrid divisions $n_m^{sg}$ = 2, 4, and 8 for a decreasing coarse cell size of $|\tilde{E}_m| \approx$ 100, 75, 50, and 25 m using the first- and second-order scheme. The error is given in a $L_2$-norm for (a) the water depth and (b) the discharge.} 
% \end{figure}

\begin{table}[htb]
    \centering
    \caption{Convergence study of a Macdonald type test with no subgrid (\textit{Classic}) and different numbers of subgrid divisions $n_m^{sg}$ = 2, 4, 8 and 16 for a decreasing coarse cell size of $|\tilde{E}_m| \approx$ 100, 75, 50, and 25 m using the first- and second-order scheme. The error is given in a $L_2$-norm for (a) the water depth and (b) the discharge.}
    \begin{tabular}{c|cccc|cccc}
         \textit{Classic} & \multicolumn{4}{c}{$P_n = 0$} & \multicolumn{4}{c}{$P_n = 1$} \\ \hline
         $\Delta s(\Omega)$ & $L_2(h)$ & $\mathcal{O}(h)$ & $L_2(hu)$ & $\mathcal{O}(hu)$ & $L_2(h)$ & $\mathcal{O}(h)$ & $L_2(hu)$ & $\mathcal{O}(hu)$ \\ \hline
         69.8  &  9.650E-02 &  -  &  2.534E-01 &  -  & 2.159E-03  &  -  &  6.431E-03  &  -    \\ 
         54.1  & 7.015E-02  & 1.3 & 1.867E-01  & 1.2 & 1.159E-03  & 2.4 &  3.418E-03  & 2.4 \\ 
         36.3  & 4.734E-02  & 1.0 & 1.275E-01  & 1.0 & 5.416E-04  & 1.9 &  1.726E-03  & 1.7 \\ 
         18.2  & 2.342E-02  & 1.0 & 6.426E-02  & 1.0 & 1.417E-04  & 1.9 &  5.029E-04  & 1.7 \\\hline  
         $n_m^{sg}=2$ & \multicolumn{4}{c}{$P_n = 0$} & \multicolumn{4}{c}{$P_n = 1$} \\ \hline
         $\Delta s(\Omega)$ & $L_2(h)$ & $\mathcal{O}(h)$ & $L_2(hu)$ & $\mathcal{O}(hu)$ & $L_2(h)$ & $\mathcal{O}(h)$ & $L_2(hu)$ & $\mathcal{O}(hu)$ \\ \hline
         69.8 &  9.530E-02 &  -  &  2.539E-01 &  -  &  2.997E-03 &  -  &  8.348E-03 &  -  \\ 
         54.1 & 6.945E-02  & 1.2 & 1.870E-01  & 1.2 & 1.851E-03  & 1.9 & 4.897E-03  & 2.1 \\ 
         36.3 & 4.698E-02  & 1.0 & 1.276E-01  & 1.0 & 1.015E-03  & 1.5 & 2.678E-03  & 1.5 \\ 
         18.2 & 2.328E-02  & 1.0 & 6.427E-02  & 1.0 & 3.453E-04  & 1.6 & 8.981E-04  & 1.6 \\ \hline
         $n_m^{sg}=4$ & \multicolumn{4}{c}{$P_n = 0$} & \multicolumn{4}{c}{$P_n = 1$} \\ \hline
         $\Delta s(\Omega)$ & $L_2(h)$ & $\mathcal{O}(h)$ & $L_2(hu)$ & $\mathcal{O}(hu)$ & $L_2(h)$ & $\mathcal{O}(h)$ & $L_2(hu)$ & $\mathcal{O}(hu)$ \\ \hline
         69.8 &  9.511E-02 &  -  &  2.545E-01 &  -  &  2.192E-03 &  -  &  7.466E-03 &  -   \\ 
         54.1 & 6.933E-02  & 1.2 & 1.873E-01  & 1.2 & 1.299E-03  & 2.1 & 4.150E-03  & 2.3 \\ 
         36.3 & 4.693E-02  & 1.0 & 1.277E-01  & 1.0 & 6.752E-04  & 1.6 & 2.220E-03  & 1.6 \\ 
         18.2 & 2.328E-02  & 1.0 & 6.430E-02  & 1.0 & 2.125E-04  & 1.7 & 6.939E-04  & 1.7 \\ \hline
         $n_m^{sg}=8$ & \multicolumn{4}{c}{$P_n = 0$} & \multicolumn{4}{c}{$P_n = 1$} \\ \hline
         $\Delta s(\Omega)$ & $L^2(h)$ & $\mathcal{O}(h)$ & $L_2(hu)$ & $\mathcal{O}(hu)$ & $L_2(h)$ & $\mathcal{O}(h)$ & $L_2(hu)$ & $\mathcal{O}(hu)$ \\ \hline
         69.8 &  9.524E-02 &  -  &  2.548E-01 &  -  &  1.945E-03 &  -  &  7.249E-03 &  -   \\ 
         54.1 & 6.942E-02  & 1.2 & 1.875E-01  & 1.2 & 1.121E-03  & 2.2 & 3.956E-03  & 2.4 \\ 
         36.3 & 4.699E-02  & 1.0 & 1.278E-01  & 1.0 & 5.519E-04  & 1.8 & 2.049E-03  & 1.6 \\ 
         18.2 & 2.330E-02  & 1.0 & 6.432E-02  & 1.0 & 1.575E-04  & 1.8 & 6.170E-04  & 1.7 \\ \hline
         $n_m^{sg}=16$ & \multicolumn{4}{c}{$P_n = 0$} & \multicolumn{4}{c}{$P_n = 1$} \\ \hline
         $\Delta s(\Omega)$ & $L_2(h)$ & $\mathcal{O}(h)$ & $L_2(hu)$ & $\mathcal{O}(hu)$ & $L_2(h)$ & $\mathcal{O}(h)$ & $L_2(hu)$ & $\mathcal{O}(hu)$ \\ \hline
         69.8 &  9.538E-02 &  -  &  2.550E-01 &  -  &  1.859E-03 &  -  &  7.187E-03 &  -   \\ 
         54.1 & 6.951E-02  & 1.2 & 1.876E-01  & 1.2 & 1.056E-03  & 2.2 & 3.896E-03  & 2.4 \\ 
         36.3 & 4.705E-02  & 1.0 & 1.279E-01  & 1.0 & 5.042E-04  & 1.9 & 1.985E-03  & 1.7 \\ 
         18.2 & 2.332E-02  & 1.0 & 6.433E-02  & 1.0 & 1.349E-04  & 1.9 & 5.917E-04  & 1.8 \\ \hline
    \end{tabular}    
    \label{tab:MacDonald_Conv}
\end{table}

The table show that the formal order of accuracy is achieved asymptotically for the first-order and second-order scheme. Moreover, the second-order scheme performs considerably better than the first-order. The first-order scheme shows no improvements when adding subgrid cells. As in this example, they get averaged out by the constant free surface and the almost constant slope in bathymetry. On the other hand, the errors are reduced for the second-order scheme when subgrid cells are added. This, we believe, is due to the approximation of the bathymetry gradient with the subgrid. Recall that the bathymetry is given on the subgrid mesh with no polynomial reconstruction, resulting in a locally first-order representation. However, the gradient is approximated more accurately with a greater number of subgrid cells. This is clear when looking at the gradient approximation using the discrete divergence theorem 
\begin{equation}
    \nabla d_m \approx \sum_{i=1}^3\sum_{(k,j)\in\mathcal{E}_{m,i}} d_{\ell(m,k)} |E_{\ell(k,m),j}| \mathbf{n}_{\ell(k,m),j}.
\end{equation}
The approximation becomes more accurate as the number of subgrid cells increases, provided the distance from the subgrid center to the coarse edge becomes smaller. This can be observed in the table, where the formal order of accuracy becomes 1.9 and 1.8 for the water depth and the discharge, respectively. In the case with $n^{sg}_m = 16$ the subgrid method is able to beat the classic second-order method for $L_2(h)$.

\subsection{Periodic solution in parabolic bowl}
\begin{figure}[htb]
    \centering
    % \vspace{3mm}
    \subfloat[]{
  \includegraphics[width=0.4\linewidth]{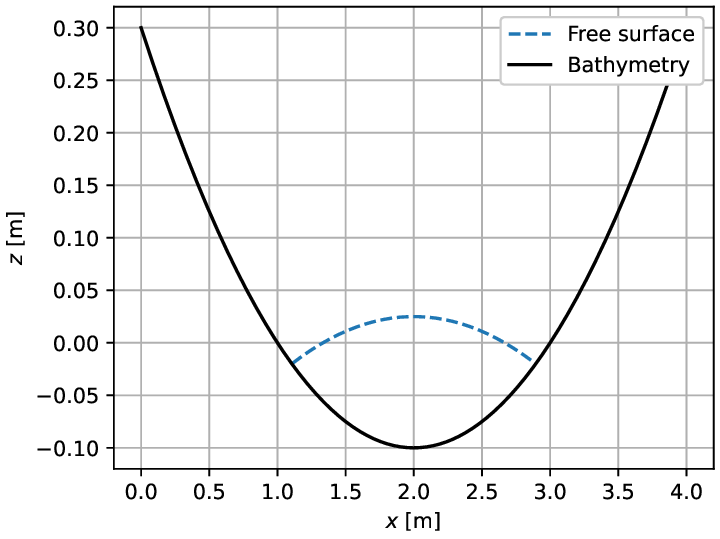} \label{fig:BM_Thacker_IC}
}
    \subfloat[]{
  \includegraphics[width=0.33\linewidth]{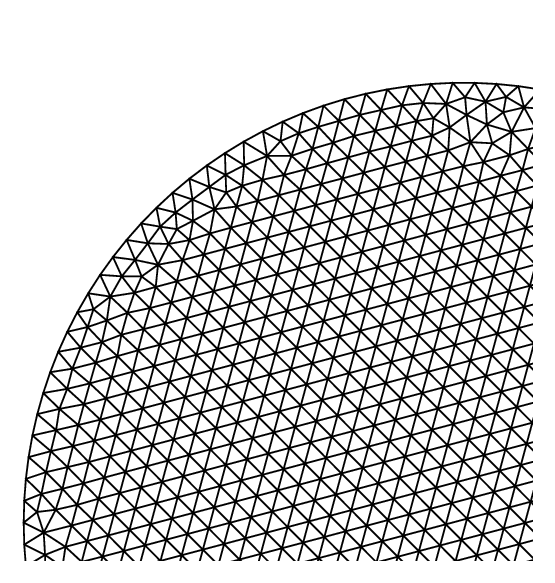} \label{fig:Thacker_Mesh}
}

    \caption{(a) Cross-sectional view at $y=0$ of the free surface initial condition and bathymetry. And (b) an illustration of one quarter of the unstructured coarse mesh used for the radially symmetric Thacker problem.}     
\end{figure}

To showcase the subgrid method's ability to accurately represent a moving flood and dry boundary, we consider a benchmark originally presented in \cite{Thacker1981}. We use the so-called radial-symmetric solution with initial condition and bathymetry showcased in Figure \ref{fig:BM_Thacker_IC}.   
% \begin{figure}[htb]
%     \centering
%     \includegraphics[width=0.55\linewidth]
%     {Figures/Thacker_Bathymetry.eps}
    
%     \caption{Cross-sectional view of the free surface initial condition and bathymetry.}   \label{fig:BM_Thacker_IC}
% \end{figure}
The water periodically runs up the bowl and down again, constantly moving the flood and dry boundary. This makes the case ideal for testing the ability to accurately capture the waterfront and demonstrate numerical damping. The case is also featured in \cite{Delestre2013} and is used in the novel subgrid work by Casulli \cite{Casulli2009}.

Due to the circular nature of the problem, a circular computational domain is used, $\Omega = [\vec{\mathbf{x}}\in \mathbb{R}^2:r=|\vec{\mathbf{x}}-(2,2)|\leq 2]$ m\textsuperscript{2}. Boundary conditions are ignorant as the water never reaches the domain boundaries. The domain is discretized using a coarse mesh consisting of 2,968 elements with an approximate element side length of $|\tilde{E}_m| \approx 1/10$ m. An illustration of one quarter of the coarse mesh is presented in Figure \ref{fig:Thacker_Mesh}.
The bathymetry is given by
\begin{equation}
    d(\vec{\mathbf{x}}) = b_0 \left( 1 - \frac{r(\vec{\mathbf{x}})^2}{a^2}\right),
\end{equation}
where $r(\vec{\mathbf{x}}) = \sqrt{(x-L/2)^2 + (y-L/2)^2}$ is the radius with center in $L/2$, and where $a$ and $b_0$ are geometrical factors. The analytical solution at any time is given by
\begin{subequations}
\begin{align}
    \eta_e(\vec{\mathbf{x}},t) &= b_0\left[\frac{\sqrt{1-A^2}}{1-A\cos{\omega t}} - 1 - \frac{r(\vec{\mathbf{x}})^2}{a^2} \left(\frac{1-A^2}{(1-A \cos \omega t)^2} -1 \right)  \right], \\
    u_e(\vec{\mathbf{x}},t) &= \frac{1}{1-A\cos{\omega t}} \left( \frac{1}{2} \omega (x - \frac{L}{2}) A \sin \omega t\right), \\
    v_e(\vec{\mathbf{x}},t) &= \frac{1}{1-A\cos{\omega t}} \left( \frac{1}{2} \omega (y - \frac{L}{2}) A \sin \omega t  \right),
\end{align}
\end{subequations}
where $\omega = \sqrt{8gb_0}/a$ is the frequency, $A = (a^2-r_0^2)/(a^2+r_0^2)$ is a constant dependent on the initial state, and $r_0$ is the distance from the center of the paraboloid to the shoreline at $t=0$. The period is given by $T = 3\frac{
2\pi}{\omega}$. The initial state is simply obtained from the analytical solution at $t=0$ s. 

\begin{figure}[htb]
    \centering
    % \vspace{3mm}
    \subfloat[$\vec{\mathbf{x}} = (2,2)$m]{
  \includegraphics[width=0.95\linewidth]{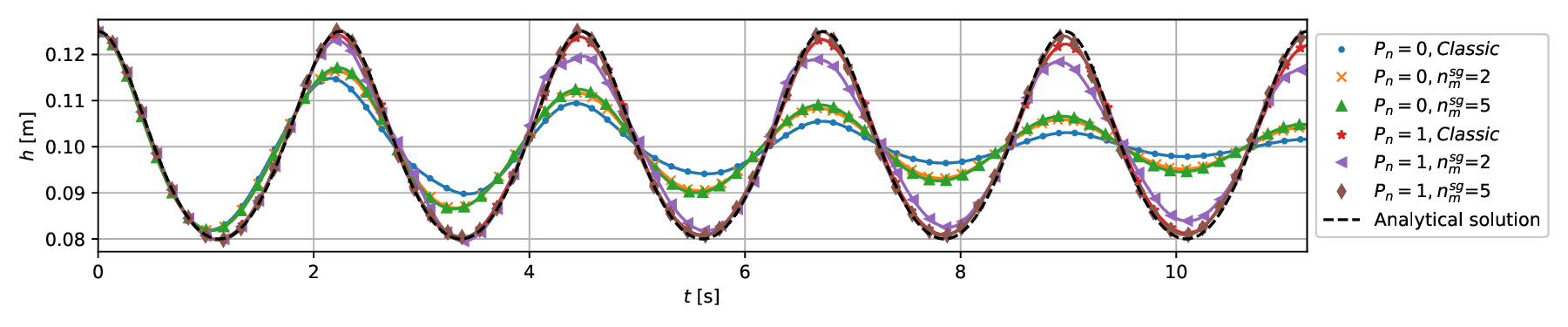} \label{fig:Thacker_2_2}
}

    \subfloat[$\vec{\mathbf{x}} = (3,2)$m]{
  \includegraphics[width=0.95\linewidth]{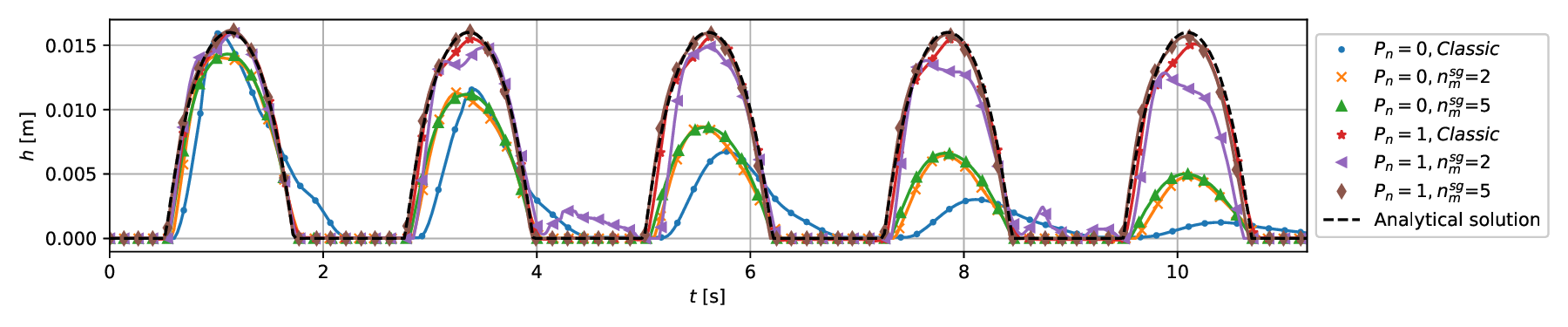} \label{fig:Thacker_3_2}
}

    \caption{Time series of the water depth recorded at two different locations for a different number of subgrid divisions and with the different polynomial order.} \label{fig:Thacker_results}
    
\end{figure}

Timeseries at $\vec{\mathbf{x}} = (2,2)$m and $\vec{\mathbf{x}} = (3,2)$ m are depicted in Figure \ref{fig:Thacker_results}, where simulations are run for $t_f = 5T$.   
The tests are run with both schemes and with no subgrid (\textit{Classic}) and subgrid divisions of $n^{sg}_m = 2$ and $5$. 
The amplitudes of the solution obtained with the first-order scheme decrease rapidly due to numerical damping. However, when comparing the cases with and without subgrid representation, the inclusion of subgrid cells leads to a significant improvement in the solution. This improvement is most evident at $\vec{\mathbf{x}} = (3,2)$m, where the domain transitions between wet and dry states, but it is also visible in the center point $\vec{\mathbf{x}} = (2,2)$m. No clear improvements are observed when increasing the number of subgrid face divisions with the first-order scheme. Similar results were reported in \cite{Casulli2009}.  
The second-order scheme performs better than the first-order scheme, as the numerical damping is reduced drastically. Both schemes are in phase at each location. The impact of subgrid is primarily visible when considering the second-order scheme. Here, the number of subgrid cells plays a vital role as $n^{sg}_m = 2$ performs worse than the true second-order scheme, but $n^{sg}_m = 5$ is more accurate compared to no subgrid (\textit{Classic}), especially for capturing the flood and dry front, see Figure \ref{fig:Thacker_3_2}. 

\subsection{Floodplain with Meandering channel}
\begin{figure}[htb]
    \centering
    \includegraphics[width=0.95\linewidth]{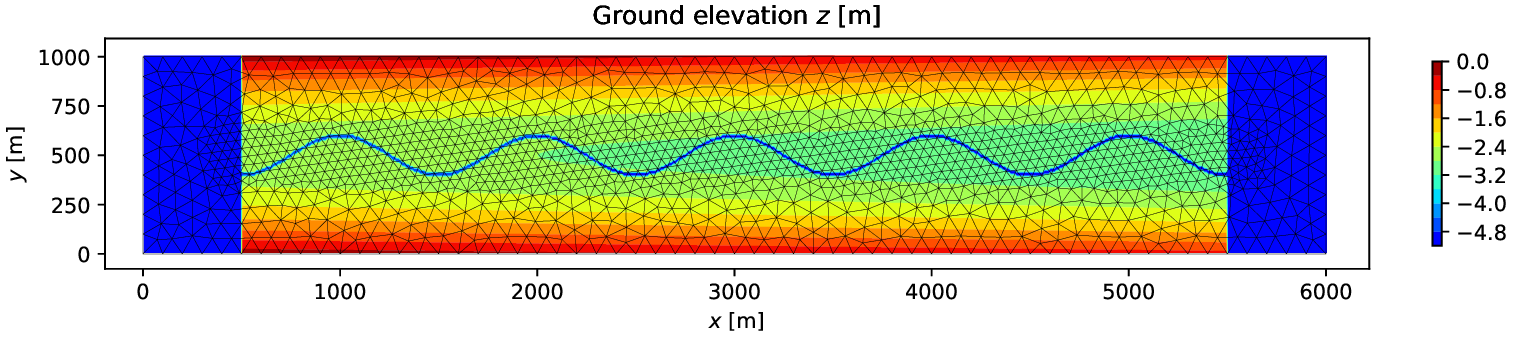} 
    \caption{Bathymetry of the parabolic floodplain with a meandering channel, and the coarse unstructured triangular computational mesh, where elements are coarse near the boundaries and finer in the middle towards the meandering channel.}
    \label{fig:PFP_Mesh}
\end{figure}

The next test is from \cite{Sanders2019} and has no analytic solution, but is designed to highlight the subgrid method's ability to resolve compound cross sections. Additionally, the test demonstrates the scheme's ability to handle discontinuous bathymetry, as the test case includes several abrupt changes in bed elevation.
The test consists of a parabolic floodplain, with two flat plateaus at each end and a small meandering channel running through the centerline of the domain (see Figure \ref{fig:PFP_Mesh}). 
The computational domain is $\Omega = [0;6000] \times [0;1000]$ m\textsuperscript{2} with a constant Manning coefficient of $M = \frac{1}{0.03}$ m\textsuperscript{1/3}/s everywhere. It is equipped with an inflow boundary condition of $\mathbf{n}\cdot h\vec{\mathbf{u}} = 2$ m\textsuperscript{2}/s at the west boundary and with wall boundary conditions elsewhere. The bathymetry of the floodplain is given as 
\begin{equation}
    z_{floodplain}(\vec{\mathbf{x}}) = -10^{-4}(x - 6000) + 10^{-5}(y-500)^2 + c,
\end{equation}
where $c$ is a constant used to shift the bathymetry for visual effects; here, we use $c = -3.2$ m, ensuring that the sign of $z_{floodplain}(\vec{\mathbf{x}})$ remains unchanged.
The depth of the meandering channel is given as
\begin{equation}
    z_{channel}(\vec{\mathbf{x}}) = -10^{-4}(x - 6000) - 2 + c  \quad 
    \text{for} \quad y \in A\sin\left((x+\frac{\lambda}{4})\frac{2\pi}{\lambda}\right) + \frac{1000}{2} \pm \frac{w_{channel}}{2},
\end{equation}
where $A=100$ m is the amplitude, $\lambda = 1000$ m is the wave length, and $w_{channel} = 20$ m is the width of the channel. The plus/minus sign should be understood as an interval such that the bathymetry should be applied when $y$ falls between the two limits. The plateaus are deeper than the rest of the domain, given by
\begin{equation}
    z_{ends}(\vec{\mathbf{x}}) = -2 + c \quad \text{for} \quad x \leq 500 \text{m  and } 5500 \text{m} \leq x.
\end{equation}
Finally, the bathymetry is given as a combination of all three parts and can be seen in Figure \ref{fig:PFP_Mesh}.
% \begin{equation}
%     d(\vec{\mathbf{x}}) = -(z_{floodplain}(\vec{\mathbf{x}}) + z_{channel}(\vec{\mathbf{x}}) + z_{ends}(\vec{\mathbf{x}}) ).
% \end{equation}

Initially, the domain is completely dry with zero velocity
\begin{equation}
    h(\vec{\mathbf{x}},0) = 0 \quad \text{and} \quad \vec{\mathbf{u}}(\vec{\mathbf{x}},0) = 0.
\end{equation}
The results are reported after $t_f = 100$ minutes.

\begin{figure}[htb]
    \centering
    \subfloat[Water depth at $t_f = 100$ minutes with $n^{sg}_m = 10$ for the first-order scheme]{
  \includegraphics[width=0.8\linewidth]{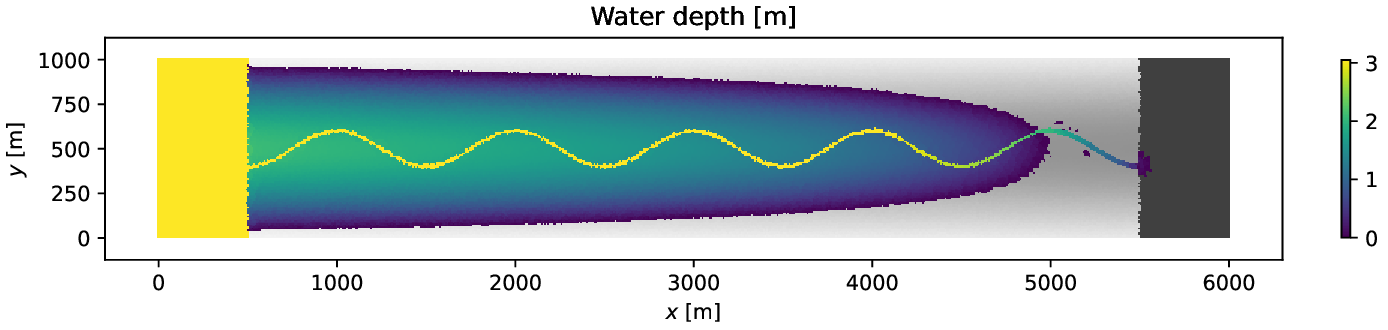} \label{fig:PFP_O1}
}

    \subfloat[Water depth at $t_f = 100$ minutes with $n^{sg}_m = 10$ for the second-order scheme]{
  \includegraphics[width=0.8\linewidth]{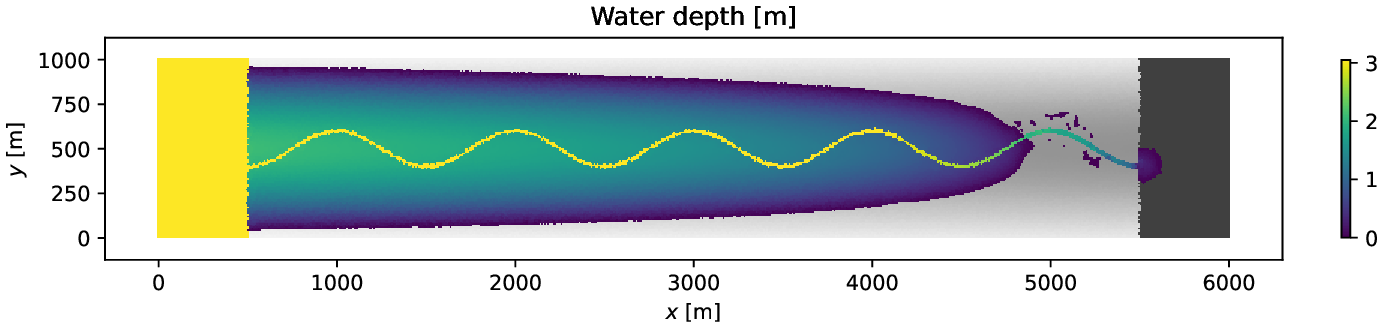} \label{fig:PFP_O2}  
}

    \subfloat[Water depth at $t_f = 100$ minutes with no subgrid for the second-order scheme]{
  \includegraphics[width=0.8\linewidth]{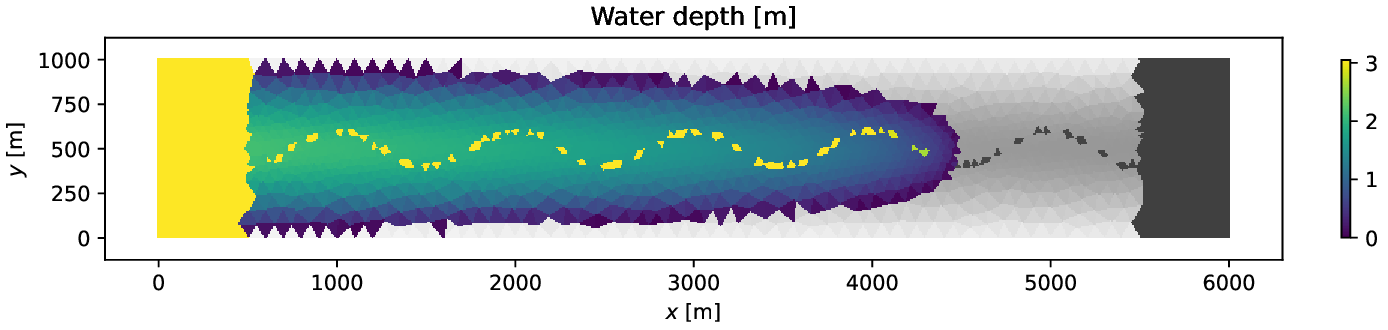} \label{fig:PFP_O2_noSG}  
}
    \caption{Top view of the water depth and bathymetry represented on the subgrid mesh, where (a) is the solution with $n^{sg}_m = 10$ and the first-order scheme, (b) is $n^{sg}_m = 10$ and with the second-order scheme, and (c) is no subgrid and with the second-order scheme.} 
\end{figure}

Here, we utilize the flexibility of the unstructured triangular mesh to have finer elements toward the meandering channel, as most of the dynamics happen here. The mesh consists of 2,876 elements, with a side length of $|\tilde{E}_m| \approx 100$ m near the boundaries and $|\tilde{E}_m| \approx 50$ m near the channel, as shown in Figure \ref{fig:PFP_Mesh}. 
Each element has a subgrid division of $n_m^{sg} = 10$, resulting in a fine subgrid element length of $|E_{\ell(m,k)}| \approx 5$ m. With this, the width of the channel is represented by 3 to 4 subgrid elements.
The results of both schemes are shown in Figures \ref{fig:PFP_O1} and \ref{fig:PFP_O2}, respectively. A solution with second-order and without subgrid is given in Figure \ref{fig:PFP_O2_noSG}.   
The test case is primarily governed by the bathymetry, and with a small variation in free surface elevation, rendering the dynamics insignificant. This is also clear when comparing the two solutions, as they are very similar. The flood and dry front in the parabolic plain is slightly behind for the second-order scheme compared to the first-order. Both the first- and second-order scheme leaks a tiny bit of water near the meandering channel. Most importantly, both results show that the water is transported at the subgrid level to the plateau at the other end, which were also observed in \cite{Sanders2019}. This point is very clear when comparing the solutions to Figure \ref{fig:PFP_O2_noSG} which is the second-order scheme without subgrid. 

\subsection{Test 7: River and floodplain linking} \label{sec:Test7}
\begin{figure}[htb]
    \centering
    % \vspace{3mm}
    \subfloat[]{
  \includegraphics[width=0.35\linewidth]{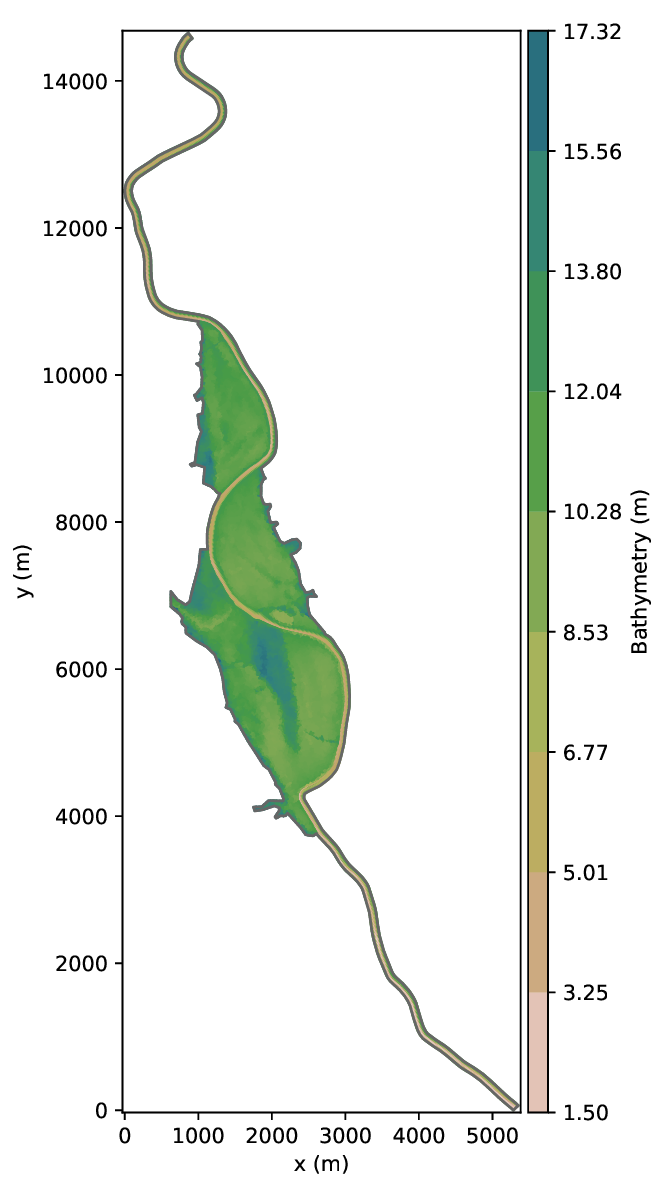} \label{fig:BM_Test7_bath}
}
    \subfloat[]{
  \includegraphics[width=0.35\linewidth]{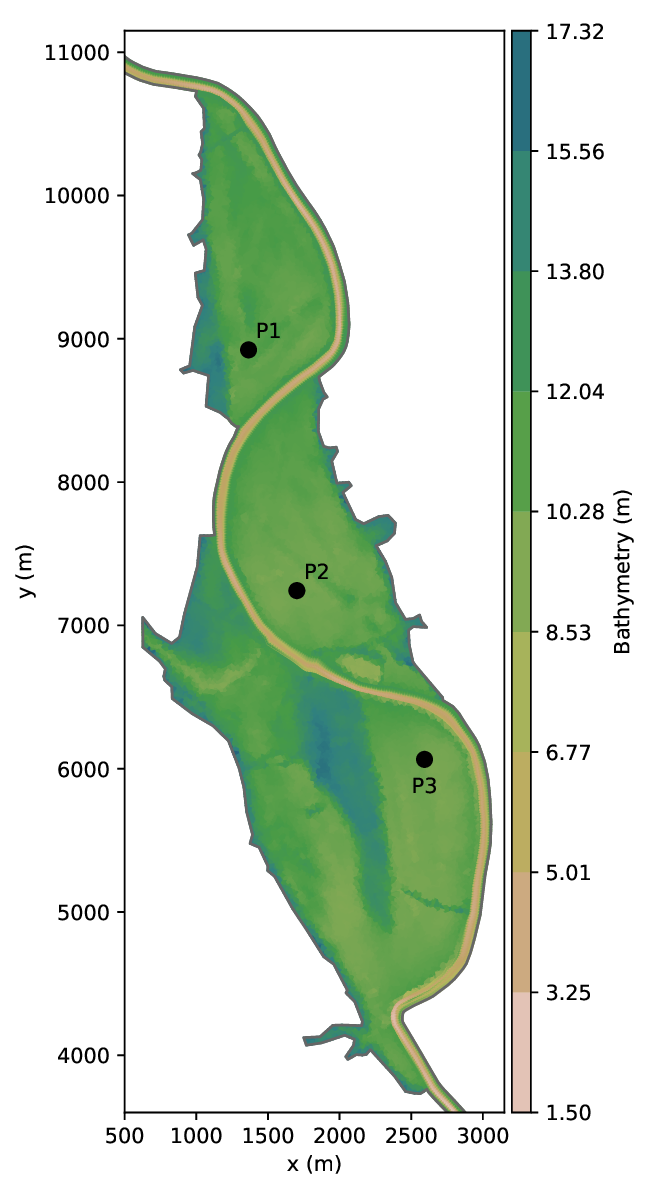} \label{fig:BM_Test7_bath_Zoom}
}

    \caption{Illustration of the computational domain: (a) the full domain, including the upstream and downstream river reaches, (b) the floodplain with the three measurement points P1, P2, and P3.}
\end{figure}

The final test case is Test 7 from the 2013 "Benchmarking the latest generation of 2D hydraulic modelling packages" by the British Environment Agency \cite{Neelz2013}. The study area covers a small section of the River Severn and the surrounding floodplain near Upton upon Severn, Worcestershire, England. In this work, a simplified version of the original test case is used, retaining the bathymetry and roughness data from the benchmark. The simplification is necessary because the subgrid model does not currently support infrastructure features such as dikes and culverts. The entire computational domain is shown in Figure \ref{fig:BM_Test7_bath}, while a zoomed-in view of the floodplains is presented in Figure \ref{fig:BM_Test7_bath_Zoom}. Furthermore, the simulations were restricted to the first-order scheme only ($P_n=0$).

The problem is initialized with a constant free surface elevation of $\eta = 9.8$ m and zero velocity, $\vec{\mathbf{u}}=\mathbf{0}$. The friction coefficient varies spatially, with $M = 35.7$ m\textsuperscript{1/3}/s in the river and $M = 25$ m\textsuperscript{1/3}/s in the floodplains.    
Flow is introduced through an upstream discharge boundary condition of 
\begin{equation}
    \mathbf{n}\cdot h\vec{\mathbf{u}}(t_{bc}) = \begin{cases}
        100 ~m^3/s, \quad t < 2 \text{ hours}, \\
        500 ~m^3/s, \quad 2 \text{ hours} \leq t < 12 \text{ hours}, \\
        100 ~m^3/s, \quad t \leq 48 \text{ hours}. \\
    \end{cases}  
\end{equation}
This setup mimics a heavy rainfall event upstream, resulting in a fivefold increase in discharge for a period of 10 hours. At the downstream boundary, water leaves the domain through a free-outflow boundary condition \cite{Sanders2002}. Simulations are performed for $t_f = 48$ hours, during which water enters the floodplains and subsequently recedes as flow conditions return towards normal.
Water depth is recorded at three locations, denoted P1, P2, and P3 (see Figure \ref{fig:BM_Test7_bath_Zoom}), corresponding to points 5, 8, and 14, respectively, in the original benchmark report.

Three meshes are used in this test case: a fine reference mesh with 49,990 elements (Mesh 1), a coarse mesh with 20,008 elements (Mesh 2), and a very coarse mesh with 6,468 elements (Mesh 3). All meshes consist of unstructured triangular elements and are constructed using two element sizes, fine elements within the river channel and coarser elements in the floodplains. The coarse meshes are illustrated in Appendix \ref{sec:appendixB}. 
The number of subgrid divisions used for each simulation is reported in Table \ref{tab:Test7}. 
For all simulations, number of time steps $N_{\Delta t}$ and the wall-clock time $t_{wall}$, was measured to demonstrate the potential computational savings by using subgrid modelling instead of conventional mesh refinement. The measured wall-clock times are reported in Table \ref{tab:Test7}.
All simulations were performed on a laptop equipped with an Intel Core Ultra 7 265H processor (16 cores, 2.2 GHz) and 64 GB of RAM. It is worth noting that all results were obtained using non-optimized Python code, and we therefore expect that the computational performance can be improved significantly through further code optimization.
\begin{table}[htb]
    \centering
    \caption{Summary of coarse meshes, subgrid configurations, number of time steps, and measured wall clock time for the different simulations in both second and hours.}
    \begin{tabular}{cccccc}
         & $N_c$ & $n^{sg}_m$& $N_{sg}$ & $N_{\Delta t} $ &$t_{wall}$ [h:m:s] \\ \hline
         Mesh 1 & 49,990 & 1 & - & 660,471 & 08:19:17 \\ \hline
     \multirow{2}{*}{Mesh 2} & \multirow{2}{*}{20,008} 
         & 1 & - & 398,652 & 02:36:46  \\ 
         & & 2 & 80,032 & 377,625 & 04:42:19  \\ \hline
     \multirow{3}{*}{Mesh 3} & \multirow{3}{*}{6,468}
        & 1 & - & 318,162 & 00:44:48 \\
        & & 2 & 25,872 & 295,476 & 01:37:11  \\
        & & 4 & 103,488 & 265,645 & 02:39:23 \\ \hline
    \end{tabular}    
    \label{tab:Test7}
\end{table}

\begin{figure}[htb]
    \centering
    \subfloat[]{
  \includegraphics[width=0.8\linewidth]{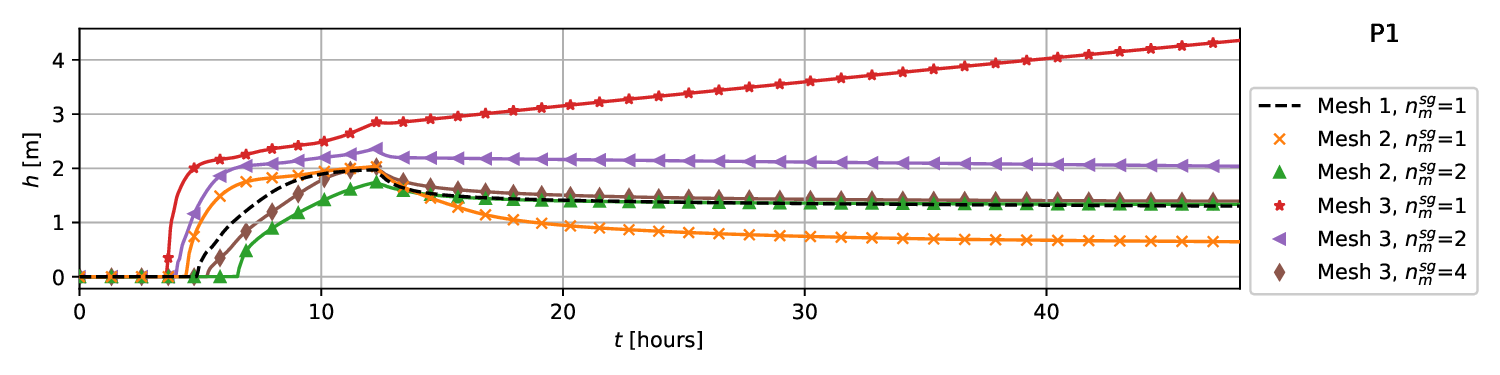} 
}

    \subfloat[]{
  \includegraphics[width=0.8\linewidth]{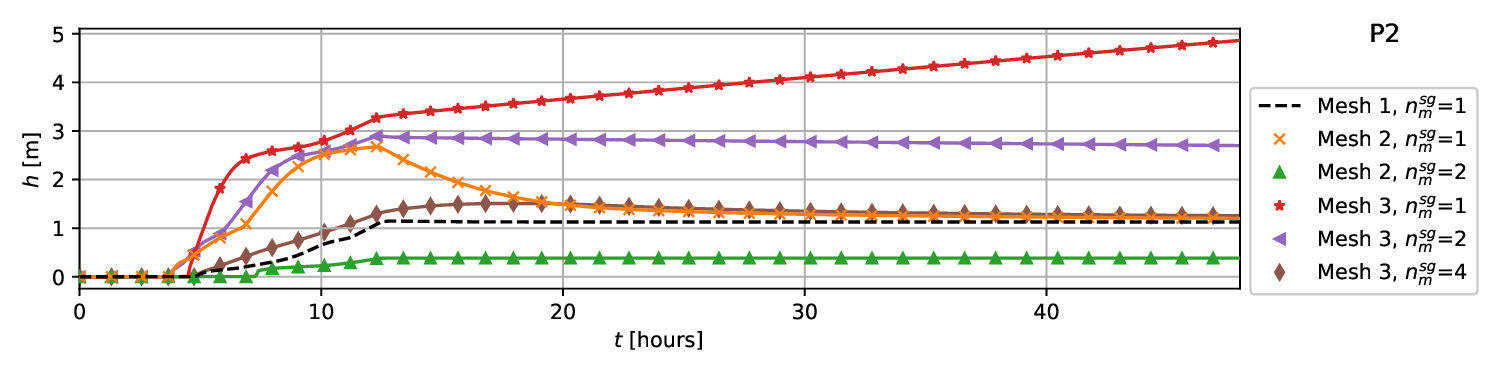} 
}

    \subfloat[]{
  \includegraphics[width=0.8\linewidth]{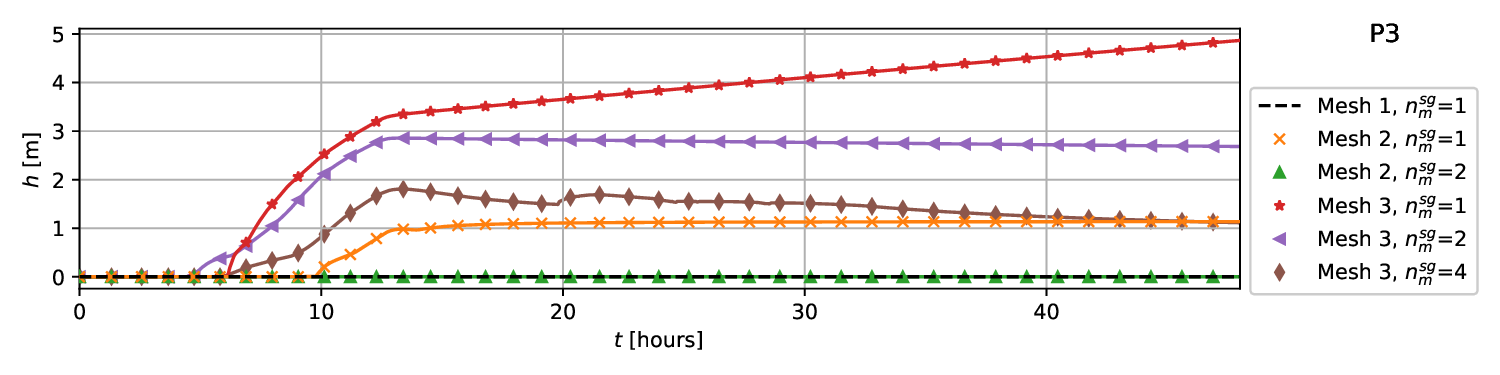}
}
    \caption{Time series of water depth at location P1 (a), P2 (b), and P3 (c). Note that the symbols are visual identifiers and not the time steps.  }\label{fig:BM_Test7_TS} 
\end{figure}

In order to compare different subgrid configurations, we consider the time series at the three measurement locations illustrated in Figure \ref{fig:BM_Test7_bath_Zoom}. The time series are shown in Figure \ref{fig:BM_Test7_TS}, where the top panel shows P1, the middle panel shows P2, and the bottom panel shows P3. Consistent with the results presented earlier, we observe that adding subgrid resolution generally improves the accuracy of the simulations, as the subgrid solutions are typically closer to the reference solution than the corresponding non-subgrid solutions. There is, however, one exception: for Mesh 2 at P2, the non-subgrid version produces a water level that is very similar to the reference solution, whereas the subgrid version lies below the reference. It is also possible that the Mesh 2 subgrid solution is, in fact, more correct than the reference solution, but investigating this possibility is beyond the scope of the present work.

\begin{figure}[htb]
    \centering
    % \vspace{3mm}
    \subfloat[Mesh 1, $n_m^{sg}=1$]{
  \includegraphics[width=0.32\linewidth]{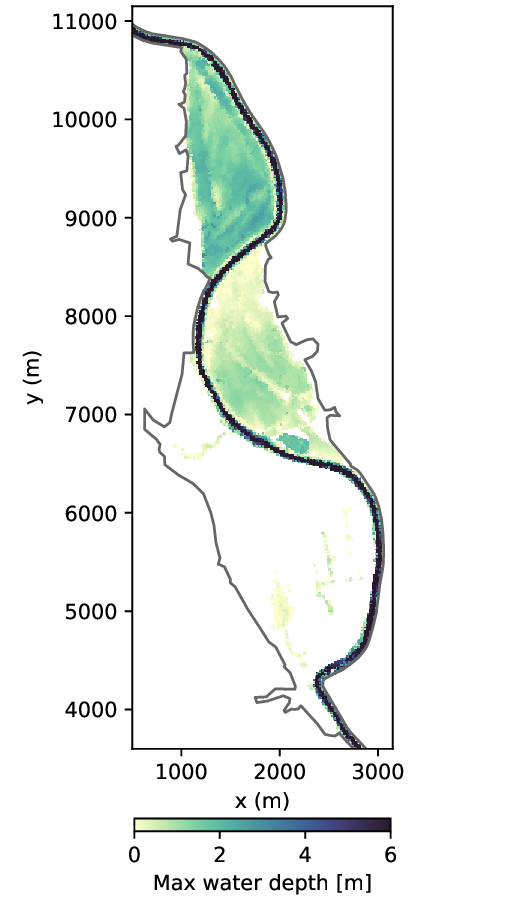} 
}
    \subfloat[Mesh 2, $n_m^{sg}=2$]{
  \includegraphics[width=0.32\linewidth]{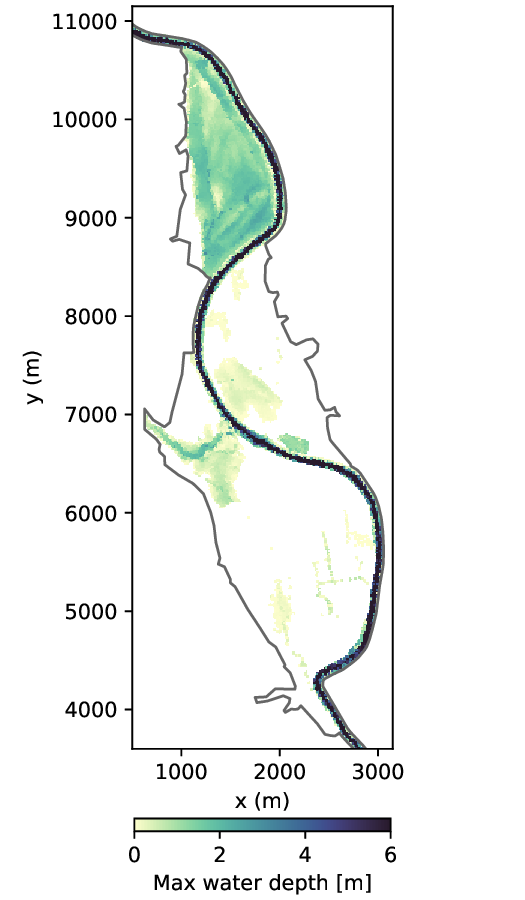} 
}
    \subfloat[Mesh 3, $n_m^{sg}=4$]{
  \includegraphics[width=0.32\linewidth]{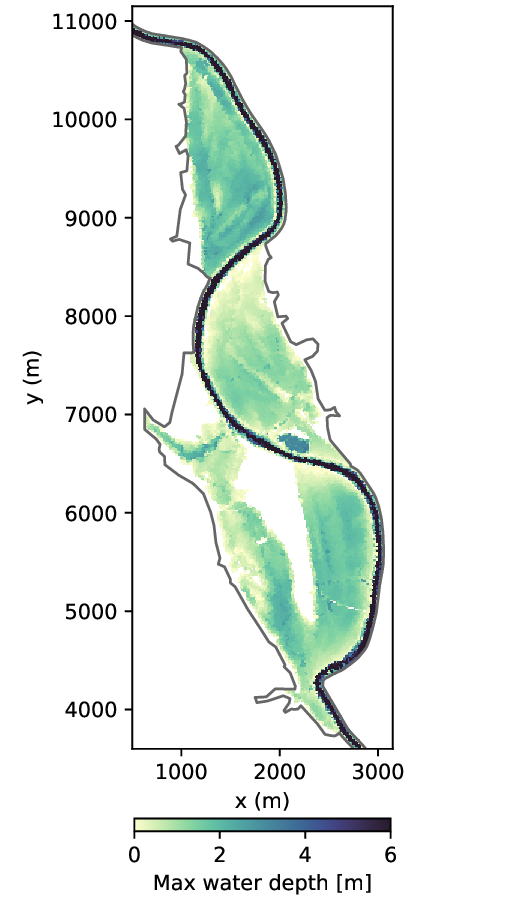}
}

    \caption{Inundation maps of the floodplains for three different mesh configurations: (a) the fine-mesh reference solution without subgrid, (b) the coarse mesh with four subgrid cells per coarse element, and (c) the very coarse mesh with sixteen subgrid cells per coarse element.} \label{fig:BM_Test7_inun}  
\end{figure}

As an additional comparison, we plotted inundation maps for three selected simulations: Mesh 1, Mesh 2 with $n_m^{sg}=2$, and Mesh 3 with $n_m^{sg}=4$, all shown in Figure \ref{fig:BM_Test7_inun}. An inundation map shows the maximum water depth reached at each location during the simulation period, thereby providing an estimate of the worst-case flooding scenario throughout the domain. Such maps are commonly used by hydraulic engineers and municipalities for flood risk assessment and urban planning. 
Comparing the two subgrid simulations with the reference solution, it is observed that the Mesh 3 simulation with $n_m^{sg}=4$ appears to overpredict inundation, as the downstream floodplain is completely flooded. However, for the upstream and middle floodplains, the results agree reasonably well with the reference solution, with only slightly greater water depths predicted in the middle floodplain. The Mesh 2 simulation with $n_m^{sg}=2$ correctly predicts that the upstream floodplain is flooded and that the downstream floodplain remains dry. However, it appears to underpredict the amount of water in the middle floodplain, at least when compared with the Mesh 1 reference solution.

Table \ref{tab:Test7} presents the recorded wall-clock times for all simulations. Despite having more subgrid elements than the reference mesh has elements, the subgrid simulations require less computational time. This is primarily due to the CFL condition being based on the coarse mesh (see Section \ref{sec:CFL}), which allows the simulations to be performed using fewer time steps, as shown in Table \ref{tab:Test7}. Furthermore, as discussed in Section \ref{sec:Model_scale}, the computational cost associated with flux calculations scales with the number of subface divisions, $n_m^{sg}$, rather than the total number of subgrid elements, $N_{sg}$.

As a conclusion to this test case we can see that we can get similar results for two different coarse meshes by adding subgrid resolution with shorter run times. However, whats important to note that all subgrid simulations start from the coarse mesh and is why it is still important to refine and calibrate this to fit the problem. This is especially true for real cases where small features can have high importance for the final results.

%% file: Texts/Conclusion_and_more.tex
\section{Conclusion} \label{sec:Conclusion}
A novel unstructured subgrid method for the SWE has been presented. The unstructured coarse mesh provides flexibility for domain fitting and local refinement of velocities. While the underlying subgrid mesh provides high bathymetry details that can capture important features that would otherwise require a very fine mesh or would be lost. 
The spatial discretization is based on the WENO method, which effectively ensures that the optimal order of accuracy is achieved. The WENO method is extended with a simple way of picking stencils near the flood and dry boundary. Together with a velocity-corrected Runge-Kutta scheme, the method is well-suited for flood and dry applications with and without friction.

The new subgrid face value reconstruction is a combination of established techniques and new subgrid cases. This takes partially wet cells into account and limits transport from coarse edges, considering both wet and dry subgrid cells. The newly developed gravity term discretization is presented and analyzed to demonstrate that it is well-balanced for both first- and second-order schemes.    

Several numerical experiments, with an emphasis on cases that have analytical solutions, have been presented to highlight the capabilities of the subgrid schemes. The cases show that key properties are preserved with subgrid, such as the order of accuracy and well-balancedness. Furthermore, it is showcased that the inclusion of subgrid can improve the solution for moving flood and dry boundaries. Finally, it is highlighted that small important features can be resolved on the subgrid level to save coarse computational mesh cells. 

Future perspectives are twofold. First, the implementation will be adapted to modern multi-core architectures, as the framework presented is parallelizable. As part of the speed-up process, we also plan to introduce a variable number of subgrid divisions based on bathymetric variation, since this is the primary driver of the proposed subgrid method. Second, we plan to achieve even higher orders of accuracy by combining the subgrid method with the Central-WENO method \cite{Dumbser2017}. To achieve a well-balanced solution, we plan to use a path-conservative scheme \cite{Pares2006,DUMBSER2009,CANESTRELLI2010}, as it has shown excellent results for high-order methods.      

\section*{Acknowledgments}
MEB acknowledges financial support from Innovation Fund Denmark (Grant no. 2052-00046). Thank you to the two anonymous reviewers for comments and questions that helped us to improve the paper.

\section*{Declarations}
The authors declare that they have no known competing financial interests or personal relationships that could have appeared to influence the work reported in this paper.

% \section*{Funding}
% The authors disclosed receipt of the following financial support for the research, authorship, and/or publication of this
% article: This work partly contributes to the activities of the research project: “Fast Flood Simulations in MIKE 21” supported by Innovationsfond (Grant no. 2052-00046).

% \begin{biog}
% \textit{Jens Visbech} is ...

% \textit{Anders Melander} is ...

% \textit{Allan Peter Engsig-Karup} is ...

% \end{biog}

%% file: Texts/Appendix.tex
\newpage

\section{Gravity source term} \label{sec:appendix}
In order to get a gravity source term that is consistent and well-balanced we start by looking at the flux and source term contributions. Recall that they are given as
\begin{subequations}
    \begin{align}
        \mathbf{F}_g(\mathbf{U}) &= \frac{1} {2}g\nabla(\eta^2 + 2\eta d), \label{eq:FluxTermG} \\
        \mathbf{S}_g(\mathbf{U}) &= g\eta \nabla d.
    \end{align}
\end{subequations}
When investigating both terms we can observe that the bathymetry gradient term is discretized in two different ways. 
To see this we start by considering the flux part
\begin{equation}
    \int_{\tilde{T}_m} \mathbf{F}_g(\mathbf{U}) ~\text{d}\vec{\mathbf{x}} = \int_{\tilde{T}_m} \frac{1}{2}g\nabla(h^2 - d^2) ~\text{d}\vec{\mathbf{x}}, \label{eq:intF}
\end{equation}
where $h^2 = \eta^2 + 2\eta d + d^2$ is used.
In the numerical scheme this is governed by the HLLC Riemann solver. Here, we have limited freedom to modify it, hence we start with some observations regarding this. In order to see the effect of the bathymetry we start by assuming that we can extract an unified water depth from the HLLC solver such that the gravity contribution can be represented as the integral in \eqref{eq:intF}. This integral can be represented with a one-point quadrature rule, making it second-order accurate, and use a discrete version of the divergence theorem to get 
\begin{equation}
    \int_{\tilde{T}_m} \mathbf{F}_g(\mathbf{U}) ~\text{d}\vec{\mathbf{x}} \approx \sum_{i=1}^3 \frac{1}{2}g((h^*_{m,i})^2 - (d^*_{m,i})^2)|\tilde{E}_{m,i}|\mathbf{n}_{m,i},
\end{equation}
where $*$ is a the unified values based on the Riemann solution, see Section \ref{sec:Flux}. 
Now, we shift the gradients to be zero at the cell center by subtracting a constant that will vanish when adding up the face contributions
\begin{equation}
\begin{aligned}
    &\sum_{i=1}^3 \frac{1}{2}g((h^*_{m,i})^2 + (d^*_{m,i})^2)|\tilde{E}_{m,i}|\mathbf{n}_{m,i} = \\ &\sum_{i=1}^3 \frac{1}{2}g(h^*_{m,i} + \bar{h}_m)(h^*_{m,i} - \bar{h}_m)| \tilde{E}_{m,i}|\mathbf{n}_{m,i} -
    \sum_{i=1}^3 \frac{1}{2}g (d^*_{m,i} + \bar{d}_m)(d^*_{m,i} - \bar{d}_m)| \tilde{E}_{m,i}|\mathbf{n}_{m,i}. \label{eq:disc_Fg}
    \end{aligned}
\end{equation}
Now to the source term; we integrate and get
\begin{equation}
    \int_{\tilde{T}_m} \mathbf{S}_g(\mathbf{U}) ~\text{d}\vec{\mathbf{x}} = \int_{\tilde{T}_m} g(h-d) \nabla d ~\text{d}\vec{\mathbf{x}},
\end{equation}
where $\eta = h - d$ is used.   
We approximate the integral with a one-point quadrature rule, the discrete divergence theorem, and shift it to get
\begin{equation}
     \int_{\tilde{T}_m} g (h-d) \nabla d ~\text{d}\vec{\mathbf{x}} \approx \sum_{i=1}^3 g \bar{h}_m(d^*_{m,i} - \bar{d}_m)|\tilde{E}_{m,i}|\mathbf{n}_{m,i} - \sum_{i=1}^3g \bar{d}_m(d^*_{m,i}-\bar{d}_m)|\tilde{E}_{m,i}|\mathbf{n}_{m,i}. \label{eq:disc_Sg}
\end{equation}
Note that we have approximated a mathematical equivalent expression 
\begin{equation}
    \frac{1}{2}g\nabla d^2 = gd\nabla d,
\end{equation}
in two different ways \eqref{eq:disc_Fg} and \eqref{eq:disc_Sg}. To clarify the bathymetry gradient term for the different expressions are given as
\begin{subequations}
    \begin{align}
        \int_{\tilde{T}_m} \frac{1}{2}g\nabla d^2 ~\text{d}\vec{\mathbf{x}} \approx& \sum_{i=1}^3 \frac{1}{2}g (d^*_{m,i} + \bar{d}_m)(d^*_{m,i} - \bar{d}_m)| \tilde{E}_{m,i}|\mathbf{n}_{m,i}, \label{eq:grad_approx1} \\
        \int_{\tilde{T}_m} g d\nabla d ~\text{d}\vec{\mathbf{x}} \approx &\sum_{i=1}^3g \bar{d}_m(d^*_{m,i}-\bar{d}_m)|\tilde{E}_{m,i}|\mathbf{n}_{m,i}.\label{eq:grad_approx2} \end{align}
\end{subequations}
This can lead to some numerical undesirable effects, when the water level is small compared to a change in the bathymetry. Here, the combined gravity can generate momentum that point uphill, which is clearly nonphysical. To avoid this we wish to approximate all gradients in the same way. Since the flux is locked due to the Riemann problem, we choose to approximate the gradients of the source term similar to the flux. This leads to the following expression for the numerical source term
\begin{equation}
     \mathcal{S}_{g}(\bar{\mathbf{U}}_m) =  \frac{1}{|\tilde{T}_m|}\frac{1}{2}g\sum_{i=1}^3 (h^*_{m,i} + \bar{h}_m)(d^*_{m,i} - \bar{d}_m)| \tilde{E}_{m,i}|\mathbf{n}_{m,i} - 
    \frac{1}{|\tilde{T}_m|}\frac{1}{2}g\sum_{i=1}^3(d^*_{m,i} + \bar{d}_m)(d^*_{m,i} - \bar{d}_m)| \tilde{E}_{m,i}|\mathbf{n}_{m,i},
\end{equation}
which in terms of of the free surface elevation is given as
\begin{equation}
    \mathcal{S}_{g}(\bar{\mathbf{U}}_m) = \frac{1}{|\tilde{T}_m|}\frac{1}{2}g\sum_{i=1}^3(\eta^*_{m,i} + \bar{\eta}_m)(d^*_{m,i} - \bar{d}_m) |\tilde{E}_{m,i}| \mathbf{n}_{m,i}.
\end{equation}
This we have shown is valid for a second-order accurate approximation in Section \ref{sec:Gravity}.   

\subsection{Unified water depth}
The gravity flux is based on the Riemann solver and face value reconstruction from the left and right side of the face. If we look at the HLLC Riemann solution one can observe that the gravity term is given as a combination of the two
\begin{equation}
    \mathcal{F}^\text{HLLC}_{g}(\bar{\mathbf{U}}_{L},\bar{\mathbf{U}}_{R},\mathbf{n}) = \alpha \mathbf{F}_g(\bar{\mathbf{U}}_R) + (1-\alpha)\mathbf{F}_g(\bar{\mathbf{U}}_L), 
\end{equation}
where $\alpha \in [0;1]$ and is based on the HLLC wave speeds
\begin{equation}
    \alpha \equiv \frac{S_R}{S_R-S_L} \quad \Rightarrow \quad (1-\alpha) = -\frac{S_L}{S_R-S_L},
\end{equation}
where $S_L$ and $S_R$ are wave speeds from the left element and right element respectively. 
From the HLLC gravity flux we get that the unified water depth is given as 
\begin{equation}
    h^2_{L,R} = \alpha h^2_R + (1-\alpha)h^2_L. 
\end{equation}

\newpage

\section{Coarse meshes for Test 7: River and floodplain linking } \label{sec:appendixB}
\begin{figure}[htb]
    \centering
    % \vspace{3mm}
    \subfloat[Mesh 1]{
  \includegraphics[width=0.30\linewidth]{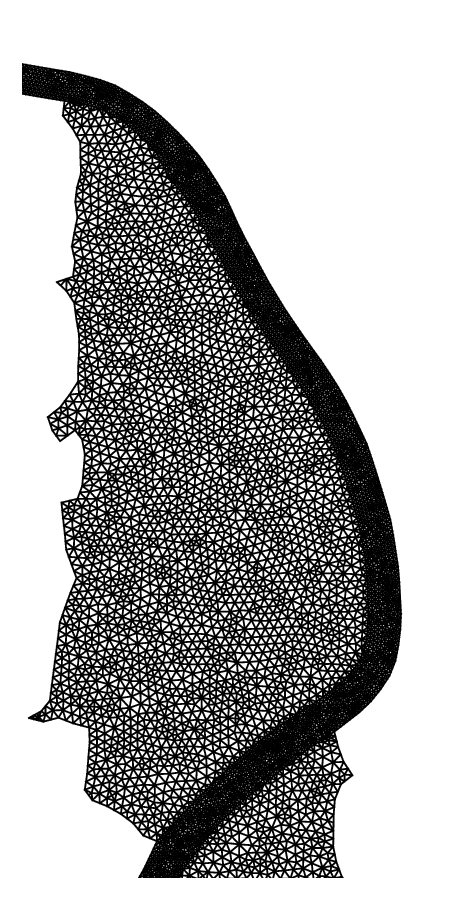} 
}
    \subfloat[Mesh 2]{
  \includegraphics[width=0.30\linewidth]{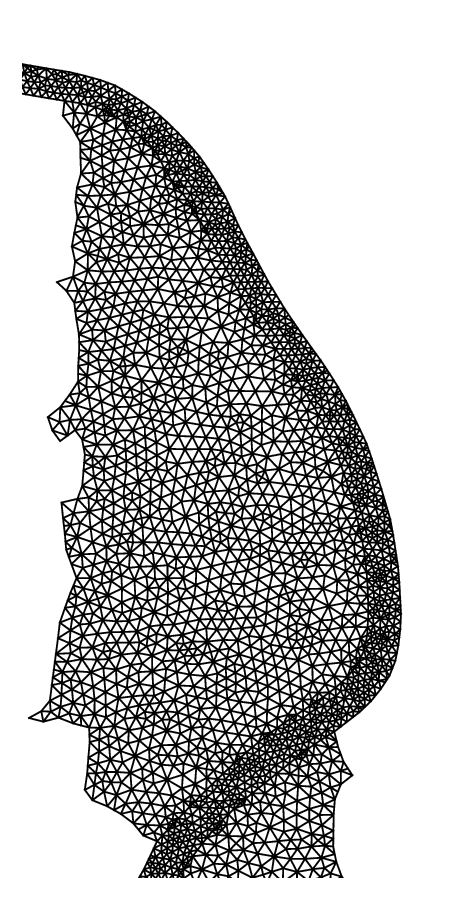} 
}
    \subfloat[Mesh 3]{
  \includegraphics[width=0.30\linewidth]{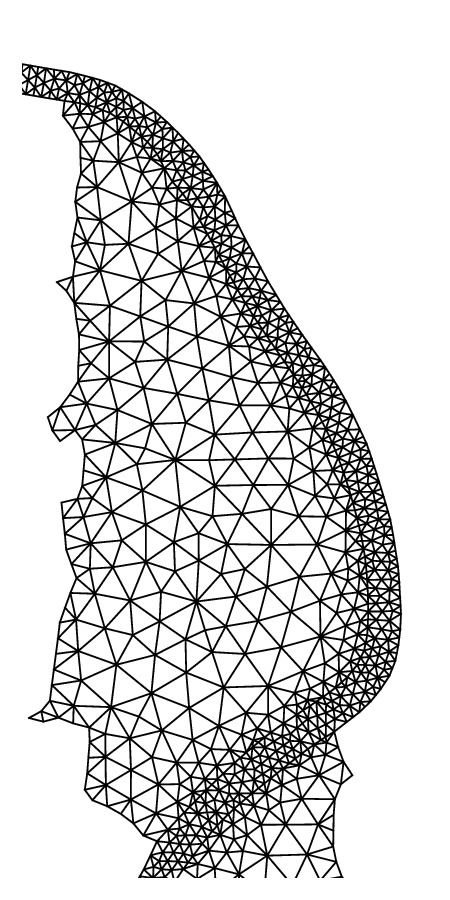} 
}

    \caption{Portions of the coarse meshes used for Test 7: River and Floodplain Linking, with (a) containing 49,990 elements, (b) 20,008 elements, and (c) 6,468 elements.} \label{fig:Test7:meshes}  
\end{figure}

Figure \ref{fig:Test7:meshes} shows the three coarse meshes used to simulate Test 7: River and Floodplain Linking (Section \ref{sec:Test7}). All meshes contain smaller elements in the river than in the floodplains to capture the higher flow velocities within the river channel. The three meshes are shown in Figure \ref{fig:Test7:meshes}, with a decreasing number of elements from left to right: 49,990 elements (Mesh 1), 20,008 elements (Mesh 2), and 6,468 elements (Mesh 3). For improved visualization, the figure is zoomed in on the uppermost floodplain. 

